\documentclass{amsart}
\usepackage{amsmath, amstext, amssymb, amsopn, amsthm, amscd, amsfonts, times,
mathptmx, epsfig}
\newcommand{\A}{\mathbb A}
\newcommand{\R}{\mathbb R}
\newcommand{\C}{\mathbb C}
\newcommand{\D}{\mathbb D}
\newcommand{\HP}{\mathbb H}
\newcommand{\M}{\mathbb M}
\newcommand{\N}{\mathbb N}
\newcommand{\Q}{\mathbb Q}
\newcommand{\Z}{\mathbb Z}

\newcommand{\T}{\mathbb T}
\newcommand{\I}{\mathbb I}

\newcommand{\SI}{\mathbb S}

\newtheorem{theo}{Theorem}

\newtheorem{prop}{Proposition}
\newtheorem{coro}{Corollary}
\newtheorem*{nntheo}{Theorem}

\newtheorem*{conj}{Conjecture}

\theoremstyle{definition}
\newtheorem{defi}{Definition}

\theoremstyle{remark}

\newtheorem{exam}{Example}

\newtheorem{note}{Note}

\newtheorem{clai}{Claim}

\newcommand{\gzg}[1]{[\negmedspace [ #1 ]\negmedspace ]}

\newcommand{\bast}{{}^{\ast}}
\newcommand{\bstar}{{}^{\text{\tiny $\bullet$}}}

\newcommand{\gcov}[1]{\gzg{\,\widetilde{ #1 }\,}}
\newcommand{\gcom}[1]{\gzg{\,\widetilde{ #1 }\,}}

\newcommand{\fgrm}{\gzg{\pi}_{1}}

\title{The Geometric Theory of the Fundamental Germ}
\author{T.M. Gendron}
\address{Instituto de Matem\'{a}ticas -- Unidad Cuernavaca, Universidad
Nacional Aut\'{o}noma de M\'{e}xico, Av. Universidad s/n, C.P. 62210
Cuernavaca, Morelos, M\'{E}XICO}
\email{tim@matcuer.unam.mx}
\date{25 May 2005}
\subjclass{Primary, 57R30, 14B0}
\keywords{fundamental germ, laminations, foliations, diophantine approximation,
mapping class group}
\begin{document}
\vspace{2cm} \maketitle
\begin{abstract} The fundamental germ 
is a generalization of $\pi_{1}$,
first defined for laminations which arise through group actions \cite{Ge1}.  In
this paper, the fundamental germ is extended to any lamination having a dense leaf
admitting a smooth structure.
In addition, an amplification of the fundamental germ
called the mother germ is constructed, which is, unlike the fundamental germ, 
a topological invariant.  The fundamental germs of the
antenna lamination and the $PSL(2,\Z )$ lamination are calculated, laminations for
which the definition in \cite{Ge1} was not available.  The mother germ
is used to give a proof of a Nielsen theorem for the
algebraic universal cover of a closed surface of hyperbolic type.
\end{abstract}

\section{Introduction}

This paper represents a continuation of our quest to extend
$\Z$-coefficient algebraic topology to laminations through the generalization
of $\pi_{1}$ called the fundamental germ.  In this paper, we extend this
construction to any lamination admitting a smooth structure.

Let us recall briefly the intuition behind 
the fundamental germ.  Consider a suspension
\[ \mathcal{L}_{\rho}\;\; =\;\;\Big( \widetilde{B}\times {\sf T}\Big) \Big/\pi_{1}B \]
of a representation $\rho :\pi_{1}B\rightarrow {\sf Homeo}({\sf
T})$, where $B$ is a manifold. 
Then $\pi_{1}B$
acts on $\mathcal{L}_{\rho}$ as fiber preserving homeomorphisms.
Let $T\approx {\sf T}$ be a fiber transversal and let $x_{0},x,\in T$.
A $\pi_{1}B$-diophantine approximation of $x\in T$ based at $x_{0}$
is a sequence $\{ g_{\alpha}\}\subset\pi_{1}B$ with 
$g_{\alpha}\cdot x_{0}\rightarrow x$.   The fundamental germ 
$\fgrm (\mathcal{L},x_{0},x)$ is then
the groupoid of tail equivalence classes of sequences of the form
$\{ g_{\alpha}\cdot h_{\alpha}^{-1}\}$ where $\{ g_{\alpha}\}$, $\{ h_{\alpha}\}$
are diophantine approximations of $x$ along $x_{0}$ \cite{Ge1}.
This construction is more generally available for any lamination 
occurring
as a quotient of a suspension, a double-coset of a Lie group or
a locally-free action of a Lie group on a space, laminations which we
refer to collectively as {\em algebraic}.
Intuitively, if $L$ is the leaf containing $x_{0}$, 
the elements of $\fgrm (\mathcal{L},x_{0},x)$ can
be thought of as sequences of paths in $L$ whose endpoints converge
transversally to $x$.  Such a sequence can be thought of as an ideal loop
based at $x$ that records an ``asymptotic identification'' within
the leaf $L$.

For a linear foliation $\mathcal{F}_{r}$ of a torus
by lines of slope $r$,
the diophantine analogy is literal and $\fgrm
(\mathcal{F},x_{0},x)$ is the group of classical diophantine
approximations of $r$.
A manifold $B$ is a supension of the trivial representation {\em i.e.}\
a lamination with a single leaf and fiber transversals that are points, which
forces $x_{0}=x$. Then all sequences in $\pi_{1}B$ converge, and we find that
$\fgrm (B,x)=\bast \pi_{1}(B,x)=$ the nonstandard fundamental group of $B$.

We now turn to the contents of this article.
The algebraic definition of the fundamental germ just described, while amenable to
calculation, has the following serious drawbacks:
\begin{enumerate}
\item  It is available only for the select family of 
algebraic laminations. 
\item  It is an invariant only with respect to the special class of
{\em trained} lamination homeomorphisms ({\em c.f.}\ \cite{Ge1}).
\end{enumerate}
Addressing these flaws is the central theme of
the present study.  In the summary that follows, we shall
assume for simplicity that all leaves are simply connected.

We begin with item (1).  
Let $\mathcal{L}$ be an arbitrary lamination 
admitting a smooth structure, let $x_{0},x$ be as above
and denote by $L$ the leaf containing $x_{0}$.  
Equip
$\mathcal{L}$ with a
leaf-wise riemannian metric that has continuous transverse variation.
In this paper, we shall refer to such a lamination as {\em riemannian}.
The new idea here is to use the leaf-wise geometry to
represent -- as sequences of isometries -- the diophantine
approximations which would make up $\fgrm$.  If 
$L$ has constant curvature geometry, this
prescription may be followed word-for-word. Fixing a transversal
$T$ containing $x_{0},x$ and a continuous section of orthonormal frames
${\sf f}=\{ {\sf f}_{y}\}$, $y\in T$, we define a diophantine
approximation of $x$ to
be a sequence $\{ A_{\alpha}\}$ of isometries of $L$ for which
$(A_{\alpha})_{\ast}{\sf f}_{x_{0}}$ belongs to ${\sf f}$
and converges transversally to ${\sf f}_{x}$.
The fundamental germ
$\fgrm
(\mathcal{L},x_{0},x , {\sf f})$ is then defined 
to be the set of tails of sequences of the form
$\big\{ A_{\alpha}B^{-1}_{\alpha}\big\}$ where
$\{ A_{\alpha}\}$, $\{ B_{\alpha}\}$
are diophantine approximations of $x$.

In the case of non constant curvature leaf-wise geometry, it is necessary
to work within the category
of {\em virtual geometry} in order to make sense of the
notion of diophantine approximation.  There, a riemannian manifold $M$ is
replaced by a union of riemannian manifolds, its virtual extension
$\bstar M$, which consists of all sequences in $M$ up to the
relation of being asymptotic. A virtual isometry between riemannian manifolds
$M$ and
$N$ consists of a pair of isometric inclusions $\bstar M
\leftrightarrows\bstar N$. All dense leaves of a riemannian
lamination have virtually isometric universal covers, and
moreover, a dense leaf having no ordinary isometries will
admit many virtual isometries.

This leads to the following definition of a diophantine
approximation: let $x\in T$, ${\sf f}$ a frame field on $T$ 
and let $L$ be any leaf accumulating on $x$.  
Then a sequence
${\sf f}_{x_{\alpha}}\rightarrow {\sf f}_{x}$, $\{ x_{\alpha}\}\subset L$, 
determines an
isometry $\bstar f:L_{x}\rightarrow U\subset \bstar
L$, where $L_{x}$ is the leaf containing $x$ and
$U$ is a component of $\bstar L$.  
The fundamental germ $\fgrm 
(\mathcal{L},L, x, {\sf f}
)$ is defined to be the set of (maximal extensions of) maps of the form
\[ \bstar f\circ \bstar g^{-1} .  \]
In this way, we now have a definition of the fundamental germ valid for
any lamination admitting a smooth structure along the leaves.

In order to address drawback (2), we will need the {\em
germ universal cover}
\[ \gcov{\mathcal{L}}\;\;\subset\;\;\bstar L ,\] defined
to be the set of asymptotic classes of sequences in $L$ that
converge to points of $\mathcal{L}$.
The germ universal cover plays the role of a unit
space for a groupoid structure on $\fgrm (\mathcal{L},L,x, {\sf f})$.
It is a lamination whose leaves are nowhere dense, and when $L$ is dense,
the canonical map $\gcov{\mathcal{L}}\rightarrow\mathcal{L}$ is onto.  
We may therefore think of
$\gcov{\mathcal{L}}$ as obtained from $\mathcal{L}$ by ``unwrapping''
all transversal topology implemented by $L$.

Assume now that $L$ is dense.  
The {\em mother germ} $\fgrm\mathcal{L}$ is defined to be the
groupoid of all partially defined maps of $\gcov{\mathcal{L}}$
that are homeomorphisms on domains which are sublaminations of $\gcov{\mathcal{L}}$ 
and 
preserve the
projection $\gcov{\mathcal{L}}\rightarrow\mathcal{L}$.  
We have in particular that
\[\fgrm \mathcal{L}\big\backslash \gcov{\mathcal{L}} \;\;\cong\;\;\mathcal{L}. \]
$\fgrm\mathcal{L}$ is
the receptacle of all the $\fgrm (\mathcal{L},L', x, {\sf f})$ for $L'$
dense, in that it contains
subgroupoids isomorphic to each.  The mother
germ is functorial with respect to topological lamination covering maps,
and is therefore, in spite of its riemannian construction, a {\em
topological} invariant.  This takes care of item (2) above.

The remainder of the paper is devoted to examples and an application.
Many examples were discussed in \cite{Ge1}, and so for this reason
we limit ourselves
to laminations which are not algebraic and hence which do not
have a fundamental germ in the sense described there.  

The first example we consider is that which we call
here the {\em antenna lamination}, a surface lamination discovered by
Kenyon and Ghys \cite{Gh}
which has the distinction of having leaves of both parabolic and hyperbolic
type.  With respect to a hyperbolic leaf, the fundamental
germ is calculated as a set to be 
$\bast F_{2}\times (\bast \Z_{\hat{2}}\oplus \bast \Z_{\hat{2}})$
where $\bast F_{2}$ is the nonstandard free group on two generators,
and $\bast \Z_{\hat{2}}$ is the subgroup of $\bast\Z$ 
isomorphic to the fundamental germ of the
dyadic solenoid.  Although a product of groups, this germ is not a group
with respect to its defined multiplication.  It is the first example we have encountered
of a fundamental
germ that is not a group.

The second example is that of the Anosov foliation of the unit tangent
bundle to the modular surface.  Although this is just the
suspension of the action of $PSL (2,\Z )$ on the boundary of the hyperbolic
plane, the definition of the fundamental germ
found in \cite{Ge1} is unavailable since it does not work for actions
with fixed points.  We calculate the fundamental germ here 
as a set to be $PSL (2,\bast \Z )$, but as in the case of the antenna lamination,
it is also not a group with respect to its defined multiplication.

The final result of this paper concerns the use of the fundamental germ
to calculate the mapping class group of the algebraic universal
cover $\widehat{\Sigma}$ of a closed surface $\Sigma$ of hyperbolic type.  
$\widehat{\Sigma}$ is by definition the inverse limit
of finite covers of $\Sigma$, a compact solenoid with dense disk leaves.
If $L\subset\widehat{\Sigma}$ is a fixed leaf, 
the leafed mapping class group ${\sf MCG}(\mathcal{L},L )$
is the quotient
${\sf Homeo_{+}}(\mathcal{L}, L )/\simeq$,
where ${\sf Homeo_{+}}(\mathcal{L}, L )$ denotes the group of
orientation-preserving homeomorphisms 
fixing set-wise $L$ and $\simeq$ denotes homotopy.
If we denote by ${\sf Vaut}(\pi_{1}\Sigma )$ the group
of virtual automorphisms of $\pi_{1}\Sigma$ ({\em c.f.}\ \S 10) then

\begin{nntheo}  There is an isomorphism
\[  \Theta :{\sf MCG}(\mathcal{L},L )\;\longrightarrow\; {\sf Vaut}(\pi_{1}\Sigma) . \]
\end{nntheo}

A proof of this theorem first appeared in the unpublished
1997 thesis
of C. Odden \cite{Od}.  Due to its importance in 
the genus-independent expression of the Ehrenpreis conjecture \cite{Ge2}, we
provide a proof in order to ensure its inclusion in the literature. 

\vspace{3mm}

\noindent
{\bf Acknowledgements:}  I would like to thank 
P. Makienko and A. Verjovsky with whom I enjoyed fruitful conversations 
regarding several important aspects of this paper.
I would also like to thank the Instituto de Matem\'{a}ticas of the UNAM
for providing a pleasant work enviroment and generous financial
support.  

\section{Virtual Geometry}

Virtual geometry is obtained as a quotient of nonstandard geometry, which
we now review: references
\cite{Gol}, \cite{Ro}, \cite{StLu}.

Let $M$ be a topological space, $\mathfrak{U}\subset {\sf 2}^{\N}$
an ultrafilter on the natural numbers all of whose elements have
infinite cardinality.  The {\em nonstandard space}
$\bast M$ is the set of sequences in $M$ modulo
$\mathfrak{U}$: that is, 
\begin{equation}\label{nonstandardspace}
\{ x_{i}\}\sim\{ y_{i}\} \text{ if and only if }\{
x_{i}\}|_{X}= \{ y_{i}\}|_{X}
\text{ for some } X\in \mathfrak{U}.
\end{equation}   
Elements of $\bast M$ are denoted $\bast x$. There is a
natural map $M\hookrightarrow \bast M$ given by the constant
sequences.  Modulo the continuum hypothesis, $\bast M$ is independent
of the choice of ultrafilter.  

There are two topologies on $\bast M$ that naturally suggest themselves.  The 
{\em enlargement topology} is generated by
sets of the form 
$\bast O$, 
where $O$ is open in $M$.  It has the
same countability as the topology of $M$ but is non-Hausdorff.  The {\em internal
topology} is generated by sets of the form
$[O_{\alpha}] =\{ \bast x\in\bast M \, | \;\bast x
\text{ is represented by a sequence }\{ x_{\alpha}\},\; x_{\alpha}\in O_{\alpha}
\}$,
where $\{ O_{\alpha}\}$ is any sequence of open sets of $M$.
It is Hausdorff but has greater countability than the topology of $M$.

For example, if we let $M=\R$ we obtain 
the {\em nonstandard reals} $\bast \R$, a totally ordered,
non-archemidean field.  Note that $\bast \R$ is
an infinite-dimensional vector space over $\R$. 
We will refer to the following substructures
of the nonstandard reals:
\begin{itemize}
\item The subring of bounded nonstandard reals, denoted
$\bast\R_{\sf fin}$, which consists of all classes of sequences
that are bounded.
\item The additive subgroup of infinitesimals, denoted
$\bast\R_{\epsilon}$, which consists of all classes of sequences
converging to $0$.
\item  The cone of positive elements, denoted $\bast\R_{+}$,
which consists of all classes of sequences that are $\geq 0$.
\end{itemize}

$\bast\R_{\sf fin}$ is a local topological ring in either the enlargement
or internal topology, with maximal ideal $\bast\R_{\epsilon}$.
The quotient $\bast\R_{\sf fin}/\bast\R_{\epsilon}$ is isomorphic to $\R$,
homeomorphic with the quotient enlargement topology (the quotient
internal topology is discrete).
The inclusion $\R\hookrightarrow\bast\R_{\sf fin}$ 
allows us to canonically identify $\bast\R_{\sf fin}$
with the product $\R\times\bast\R_{\epsilon}$.  Taking
the product of the euclidean topology on $\R$ with
the discrete topology on $\bast\R_{\epsilon}$,
we obtain a third topology on $\bast\R_{\sf fin}$ which
is Hausdorff and
quotients by $\bast\R_{\epsilon}$ to the topology on $\R$.
We call this third topology the {\em lamination topology}:
it may be extended to $\bast \R$ by giving the group $\bast\R/\R$
the discrete topology and identifying $\bast\R\cong\R\times (\bast\R/\R )$.

If $M$ is an $n$-manifold, then $\bast M$ is a nonstandard manifold
modelled on $\bast\R^{n}$.
If we denote by $\bast M_{\sf fin}$ the points of $\bast M$ represented
by sequences which converge 
to points of $M$, then we may choose an atlas on $\bast M_{\sf fin}$ whose
transitions preserve the lamination structure of $\bast \R^{n}_{\sf fin}$ {\em i.e.}\
$\bast M_{\sf fin}$ is an $n$-lamination.  In general, $\bast M$ is
a union of laminations of dimensions $\leq n$, 
this because of the possibility of ``dimension collapse''
which we describe in the proof of Theorem~\ref{virtualmanifold} below.

If $d$ is a metric inducing the topology of $M$, it extends to a
$\bast\R_{+}$-valued metric $\bast d$ on $\bast M$. Write $\bast
x\simeq\bast x'$ if $\bast d(\bast x ,\bast x' )\in \bast
\R_{\epsilon}$.

\begin{defi}  The {\bf virtual extension} of $M$ is
the quotient
\[      \bstar M\;\; =\;\; \bast M/\simeq \, ,       \]
equipped with the quotient lamination topology.  
\end{defi}

The virtual extension of $\bstar\R$ of $\R$ is called the {\em virtual reals},
a totally-ordered real vector space.
The metric $\bast d$ on $\bast M$ induces a $\bstar\R_{+}$-valued
metric $\bstar d$ on $\bstar M$.  Given $\bstar x\in\bstar M$, the
set 
\[  U_{\bstar x}\;\; =\;\; \{ \bstar y\; |
\;\; \bstar d(\bstar x ,\bstar y)\in\R  \}  \]
is a component of $\bstar M$ called the 
{\em galaxy} of $\bstar x$.  $M$ is a galaxy of $\bstar M$, and $\bstar M$ is
the union of all of its galaxies.

The galaxies of $\bstar M$ can be quite different
from one another.  For example
if $M$ is simply connected, then there may be galaxies that are not.
For example, suppose that $M$ is a noncompact leaf of the Reeb foliation
of the torus.  Consider a sequence of points $\{ x_{\alpha}\}$ in 
$M$ converging to a point $\hat{x}$ in the compact toral leaf.  
Let  $\{ \gamma_{\alpha} \}$ be a sequence of simple closed curves converging
to the meridian through $\hat{x}$.  Then the limit curve $\bstar \gamma$
is essential in the universe $U_{\bstar x}$.
On the other hand, if $M$ is a riemannian homogeneous
space, then the universes of $\bstar M$ are all isometric to $M$.

\begin{theo}\label{virtualmanifold}  If $M$ is a complete riemannian manifold
of dimension $n$, each galaxy $U$
of $\bstar M$ has the structure of 
a complete riemannian manifold of dimension $m\leq n$.
\end{theo}

\begin{proof}  Given a galaxy $U$, $\bstar x\in U$ and $\{ x_{\alpha}\}$
a representative sequence, let $m$ be the largest integer
for which there exists a sequence of $m$-dimensional balls
$\{ D_{r}(x_{\alpha})\}$ of fixed radius $r$ about $\{ x_{\alpha}\}$.  
The integer $m$ is independent of the representative sequence 
and defines an $m$-ball $D_{r}(\bstar x )\subset U$.  The
function $\bstar x\mapsto m$ is locally constant, thus 
the collection of such balls defines on  
$U$ the structure of a smooth $m$-manifold. 
Note that it is possible to have $m<n$: for example, if 
$M$ is a hyperbolic manifold with a cusp, then for
a class of sequence emptying into the cusp, we have 
$m=n-1$. 

Consider the nonstandard tangent bundle
\[ {\bf T}\bast M\;\; :=\;\; \bast\big({\bf T}M\big). \]
There is a natural projection of ${\bf T}\bast M$ onto $\bast M$
whose fiber ${\bf T}_{\bast x}\bast M$ -- the tangent space at $\bast x$ --
consists of classes of sequences of vectors $\{{\sf v}_{\alpha}\}$
based at sequences $\{ x_{\alpha}\}$ belonging to the class of $\bast x$.
It is not difficult to see that ${\bf T}_{\bast x}\bast M$ is a real 
infinite-dimensional vector space.  The riemannian
metric $\rho$ extends to a $\bast\R$-valued metric $\bast\rho$ on 
${\bf T}\bast M$ in the obvious way.  Denote by $\bast |\cdot |$ the associated
norm.  Define the bounded tangent bundle by
\[  {\bf T}_{\sf fin}\bast M\;\; =\;\;
\Big\{ \bast {\sf v}\in {\bf T}\bast M \;  
\Big| \;\; \bast | \bast {\sf v}|\in \bast\R_{\sf fin}  \Big\} . \]

Given tangent vectors $\bast {\sf v}$ and $\bast {\sf v}'$ based
at $\bast x$ and $\bast x'$, 
we write $\bast {\sf v}\simeq\bast {\sf v}'$ if 
\begin{enumerate}
\item $\bast x\simeq\bast x'$. 
\item the Levi-Civita parallel translate of a representative
$\{{\sf v}_{\alpha}\}$ of $\bast {\sf v}$ to a representative $\{ x_{\alpha}'\}$ 
of $\bast x'$ -- 
along a sequence of geodesics connecting to a representative 
$\{ x_{\alpha}\}$ of $\bast x$ -- is asymptotic to a representative
$\{{\sf v}_{\alpha}'\}$ of $\bast {\sf v}'$.
\end{enumerate}

Now define the bounded tangent bundle of $\bstar M$ to be
\[  {\bf T}_{\sf fin}\bstar M\;\;=\;\; {\bf T}_{\sf fin}\bast M\big/ \simeq .\]
The nonstandard riemannian metric $\bast\rho$ on ${\bf T}_{\sf fin}\bast M$
descends to a riemannian metric on ${\bf T}_{\sf fin}\bstar M$.
If $U$ is a galaxy, its tangent space  
may be identified with the restriction of ${\bf T}_{\sf fin}\bstar M$ to
$U$.  
Now any geodesic $\eta\subset U$ can be realized as a sequence class of geodesics
$\{ \eta_{\alpha}\}$.  Since each member of such a sequence can be continued 
indefinitely, the same is true of $\eta$, hence $U$ is complete.
\end{proof}

\begin{defi}  Let $M$, $N$ be riemannian $n$-manifolds.  A {\bf virtual subisometry}
is an injective map
\[ \bstar f:\bstar M\hookrightarrow \bstar N,  \]
where $\bstar f$ maps each galaxy of $\bstar M$ isometrically
onto a galaxy of $\bstar N$.  If in addition there exists a
virtual subisometry $\bstar g :\bstar N\hookrightarrow \bstar M$,
then the pair $(\bstar f, \bstar g )$ is called a {\bf virtual isometry}.
\end{defi}

We write $M\leq_{\sf vir} N$ to
indicate the existance of a virtual subisometry $\bstar f$ and
$M\cong_{\sf vir}N$ indicates the existence of a virtual isometry.  The relation
$\leq_{\sf vir}$ defines a partial ordering on the set of all
riemannian $n$-manifolds.

An isometry $f:M\rightarrow N$ clearly induces a virtual isometry
$(\bstar f ,\bstar g ):\bstar M\leftrightarrows \bstar N$ with 
$\bstar f ,\bstar g $ inverse to one
another.  More generally, a continuous map $\bstar f:\bstar M\rightarrow\bstar N$ 
is called {\em standard} if it is induced by a map $f:M\rightarrow N$ {\em i.e.}\
if for any $\bstar x\in \bstar M$ and any representative $\{ x_{\alpha}\}$, $\{ f(x_{\alpha})\}$
is a representative of $\bstar f(\bstar x )$.

\begin{theo}\label{leavesvirtrelated}  Let $L$ be a dense leaf of a riemannian 
lamination $\mathcal{L}$.  Then for every leaf $L'\subset\mathcal{L}$, 
\[ \widetilde{L}'\;\leq_{\sf vir}\;\widetilde{L}.\]
\end{theo}

\begin{proof}  Fix a global metric $d$ on $\mathcal{L}$
which agrees locally with the riemannian metric on the leaves.
(By this we mean that in sufficiently small flow boxes,
$d$ agrees with the distance function of $\rho$
in any plaque.)  Let $\{ \tilde{x}_{\alpha}'\}\subset \widetilde{L}'$ be any sequence,
$\{ x_{\alpha}'\}$ its projection to $L'$.
Let $\{ \tilde{x}_{\alpha}\}\subset \widetilde{L}$ be a sequence whose
projection $\{ x_{\alpha}\}$ to $L$ is $d$-asymptotic to 
$\{ x_{\alpha}'\}$. By transversal
continuity of the metric, we deduce a sequence of
$K_{\alpha}$-quasiisometries, $K_{\alpha}\rightarrow 1$,
\[ f_{\alpha}:D_{\delta}(x_{\alpha})\rightarrow D_{\delta}(x_{\alpha}'),\]
for some $\delta >0$, where $D_{\delta}(x )$ means the open 
$\rho$-ball of radius $\delta$ about $x$.  Then if $\bstar \tilde{x}$, $\bstar \tilde{x}'$
are the virtual classes of $\{ \tilde{x}_{\alpha}\}$, $\{
\tilde{x}'_{\alpha}\}$, the sequence of quasiisometries $\{ f_{\alpha}\}$ 
induces an isometry $D_{\delta}(\bstar
\tilde{x})\rightarrow D_{\delta}(\bstar \tilde{x}')$. Since $L$ is dense, we may
continue these isometries along geodesics to obtain a locally isometric surjection
$U\rightarrow  U'$, where $U$, $U'$ are the galaxies containing
$\bstar \tilde{x}$, $\bstar \tilde{x}'$.  But since these spaces
are simply connected, and the map is isometric, this surjection is a bijection.  
Hence it inverts to an isometry $U'\rightarrow U$.
Repeating this for every
$\bstar$-class of sequence in $\widetilde{L}'$,
we obtain the desired virtual subisometry 
$\widetilde{L}'\leq_{\sf vir}\;\widetilde{L}$.
\end{proof}

Two riemannian manifolds have the same {\em virtual geometry} if their universal covers
are virtually isometric.

\begin{coro}  Dense leaves of a riemannian lamination $\mathcal{L}$ have the same
virtual geometry.
\end{coro}

\section{The Fundamental Germ}

Let $\mathcal{L}$ be a riemannian lamination, $x$ a point contained
in a transversal $T$, $L$ a leaf accumulating at $x$ and $L_{x}$
the leaf containing $x$.
Let ${\sf f} :T\rightarrow {\bf
F}_{\ast}\mathcal{L}$ be a continuous section of the leaf-wise orthonormal frame bundle
of $\mathcal{L}$
over $T$.  Fix locally isometric universal
covers $p:\widetilde{L}\rightarrow L$ and $p_{x}:\widetilde{L}_{x}\rightarrow L_{x}$. 
Denote $T_{0}=T\cap L$, $\widetilde{T}_{0}=p^{-1}(T_{0})$ and let
$\tilde{\sf f}_{\tilde{y}}$ denote the lift of ${\sf f}_{y}$ to 
a point $\tilde{y}\in\widetilde{T}_{0}$ covering $y$.  We pick a basepoint
$\tilde{x}\in\widetilde{L}_{x}$ lying over $x$ with lifted frame
$\tilde{\sf f}_{\tilde{x}}$.

Let $\tilde{y}\in \widetilde{T}_{0}$.  For $r>0$, the frames 
$\tilde{\sf f}_{\tilde{x}}$, $\tilde{\sf f}_{\tilde{y}}$ determine
polar coordinates on the
metric disks $D_{r}(\tilde{x})$, $D_{r}(\tilde{y})$.  This yields
in turn a canonical quasiisometry 
\[ f: D_{r}(\tilde{x})\rightarrow D_{r}(\tilde{y})\]
given by the coordinate maps.

Let $\{ x_{\alpha}\}\subset T_{0}$ be a sequence converging to $x$, 
$\{\tilde{x}_{\alpha}\}\subset \widetilde{T}_{0}$ any sequence covering 
$\{ x_{\alpha}\}$. 
Then the frame sequence $\{ \tilde{\sf f}_{\tilde{x}_{\alpha}}\}$ and the frame
$\tilde{f}_{\tilde{x}}$ determine 
a sequence of $K_{\alpha}$-quasiisometries 
\[ \Big\{ f_{\alpha}
:D_{r_{\alpha}}(\tilde{x})\;\longrightarrow\; 
D_{r_{\alpha}}(\tilde{x}_{\alpha})\Big\} .\] 
Since $L$ accumulates at $x$, we may choose the sequence of
radii $r_{\alpha}\rightarrow\infty$
so that $K_{\alpha}\rightarrow 1$.
We deduce an isometry
\[  \bstar f:  \widetilde{L}_{x}\longrightarrow U\subset\bstar\widetilde{L} \]
where $U$ is the galaxy containing $\bstar \tilde{x}$. 
The map $\bstar f$ 
is called an {\em {\sf f}-diophantine approximation} of $x$
along $L$.

\begin{defi}\label{fundgerm}  The {\bf fundamental germ} of $\mathcal{L}$, 
based at $x$ along $L$ and ${\sf f}$, is
\[ \fgrm (\mathcal{L},L,x,{\sf f} )\;\; =\;\; \Big\{ \bstar f\circ
\bstar g^{-1}\;\Big|\;\; \bstar f,\, \bstar g\text{ are } {\sf
f}\text{-diophantine approximations of }x\text{ along } L \Big\}. \]
\end{defi}

If $x\in L$, we shorten the notation to $\fgrm (\mathcal{L}, x,{\sf f} )$.
The groupoid structure of $\fgrm (\mathcal{L},L,x,{\sf f} )$ will be described
in the next section.

\begin{note}  Suppose that $\mathcal{L}$ is a constant curvature
riemannian foliation
with dense leaf $L$ modeled on the space form $\M^{n}=\R^{n}$ or $\HP^{n}$.
Then the frame field actually determines
a sequence of uniquely defined {\em global} isometries $\{ f_{\alpha}:
\M^{n}\rightarrow \M^{n}\}$.  
Given $G$ a group, nonstandard $G$ is the group $\bast G$ of all
sequences $\{ g_{\alpha}\}\subset G$ modulo the relation $\sim$ described
in (\ref{nonstandardspace}).
Then an ${\sf f}$-diophantine approximation is completely determined
by the class $\bast f\in\bast {\sf Isom}(\M^{n})$ of $\{ f_{\alpha}\}$.
We note that $\bast {\sf Isom}(\M^{n})$
is a subgroup of $ {\sf Isom}(\bstar \M^{n} )$ 
(the group of isometries of $\bstar\M^{n}$, {\em not}
virtual isometries).  Thus, if $\Gamma \cong \pi_{1}L$ is the deck
group of $\M^{n}\rightarrow L$, we have
\[  \bast\Gamma \;\; \subset\; \;  \fgrm (\mathcal{L},L,x, {\sf f} )
\;\;\subset\;\;    \bast {\sf Isom}(\M^{n} ) .\] 
\end{note}

The terminology 
${\sf f}$-diophantine approximation comes from the following example.

\begin{exam}\label{diophex}   Let
$\mathcal{L}$ be the irrational foliation of the torus $\T^{2}$ by
lines of slope $r\in\R\setminus\Q$.  Define a representation 
$\rho :\Z \cong \pi_{1} \SI^{1}\rightarrow{\sf Homeo}(\SI^{1})$ by
$\rho_{n}(\bar{y})=\overline{y-nr}$, where $\bar{y}$ denotes the image
of $y\in\R$ in $\SI^{1}=\R/\Z$.  Then the suspension of $\rho$,
$\mathcal{L}_{\rho}= (\R\times \SI^{1})/\Z$ , is
homeomorphic to $\mathcal{L}$.  The map $\R\times\SI^{1}\rightarrow\SI^{1}$
defined $(x,\bar{y})\mapsto\bar{x}$ ({\em i.e.}\ the projection onto
the first factor composed with the universal
covering $\R\rightarrow\SI^{1}$) induces
a projection 
$\mathcal{L}_{\rho}\rightarrow \SI^{1}$.  Let
$T\approx \SI^{1}$ be a fiber of this
projection passing through
$x$.  A frame section ${\sf f}$ along $T$ is determined by an
orientation of $\mathcal{L}$. In this case, an ${\sf
f}$-diophantine approximation of $x$ is just a diophantine
approximation of $r$.  (Recall that a sequence $\{ n_{\alpha}\}\subset\Z$ is called a
diophantine approximation of $r\in\R$ if $\{ \overline{rn_{\alpha}}\}$
converges to $\bar{0}\in\SI^{1}$.)  Thus if one denotes by 
$\bast\Z_{r}$ the subgroup of $\bast\Z$ consisting of classes of
diophantine approximations of $r$, we
obtain in agreement with the construction in \cite{Ge1}, \S 4.4:
\[ \fgrm (\mathcal{L},L,x,{\sf f} )\;\; =\;\; \bast\Z_{r}. \]
(Note: $\bast\Z_{r}$ is an ideal if and only if $r$ is rational.)
If another frame
field ${\sf f}'$ is used whose domain is a transversal $T'$ which is not a
suspension fiber, the set of diophantine approximations is a subset
 $\bast\R_{r}\subset\bast\R$.  This subset maps injectively
into $\bstar\R$ with image $\bstar\Z_{r}=\bast \Z_{r}$.
\end{exam}

\begin{exam}\label{invlimitsol}  Consider a nested set of Fuchsian groups 
$\mathcal{G}=\{ \Gamma_{i}\}$
and let
\[ \widehat{\Sigma}_{\mathcal{G}}\;\; =\;\;\lim_{\longleftarrow}\;\HP^{2}/\Gamma_{i} ,\]
be the associated hyperbolic surface solenoid.  We may take $T$ to
be a fiber $\hat{p}^{-1}(x_{0})$ of the projection
$\hat{p}:\widehat{\Sigma}_{\mathcal{G}}\rightarrow \Sigma_{0}$, where
$\Sigma_{0}=\HP^{2}/\Gamma_{0}$ is the initial
surface. Then a frame at $x_{0}$ pulls back to a frame section
${\sf f}$ along $T$.  In this case, we find that
Definition~\ref{fundgerm} again agrees with the definition
found in \cite{Ge1}:
\begin{eqnarray*}
\fgrm (\widehat{\Sigma}_{\mathcal{G}},L,x,{\sf f} )\;\; & = &\;\;
\bigcap\bast\Gamma_{i} \\
& = & \;\;
\Big\{ \{ g_{\alpha}\}\subset\Gamma_{0}\;\Big|\;\; \text{for all } i,\;\exists\;
N_{i}\text{ such that }
g_{\alpha}\in \Gamma_{i}\text{ when }\alpha>N_{i}\Big\}\Big/\sim,
\end{eqnarray*}
a subgroup of $\bast PSL (2,\R )\cong PSL( 2,\bast\R )$. If ${\sf f}'$ is another
frame field, not necessarily with a fiber transversal domain, then
the corresponding germ $\fgrm
(\widehat{\Sigma}_{\mathcal{G}},L,x,{\sf f}' )$ need not define a subgroup
of $PSL( 2,\bast\R )$
and particularly, need not be
isomorphic to
$\bigcap\bast\Gamma_{i}$ (although the fundamental germs
calculated with respect to ${\sf f}$ and ${\sf f}'$ are in
canonical bijection). The rub here is the non-uniform nature of
the action of $PSL(2,\R )$ on $\HP^{2}$.  This problem will become
moot through the replacement of the fundamental germ by the mother germ, 
\S\ref{mother}.
\end{exam}

\begin{exam}  More generally, let $\mathcal{L}$ be any hyperbolic surface
lamination.  Then the fundamental germ
$\fgrm (\mathcal{L},L,x, {\sf f} )$ is a subset of
$PSL(2,\bast\R )$. Equally, if $\mathcal{L}$ is a hyperbolic 3-lamination,
$\fgrm (\mathcal{L},L,x, {\sf f} )\subset
PSL(2,\bast \C )$.
\end{exam}

\section{The Germ Universal Cover}\label{germuniversalcover}

Let
$L\subset \mathcal{L}$ be a fixed leaf.  Denote by
$p:\widetilde{L}\rightarrow L$ the universal cover.  We recall the following
definition \cite{Ge1}:

\begin{defi}  
The {\bf germ universal cover} of $\mathcal{L}$ along $L$
is the subspace $\gcov{\mathcal{L}}\subset\bstar\widetilde{L}$ defined 
\[ \gcov{\mathcal{L}}
\;\;=\;\;\Big\{ \{ \tilde{x}_{\alpha}\}\subset \widetilde{L}\;  \Big|\;\;
\big\{ p(\tilde{x}_{\alpha})\big\}\text{ converges in }\mathcal{L}
\Big\}\Big/ \simeq . \]
\end{defi}
We will denote elements of the germ universal cover by $\bstar \tilde{x}$.
There is a natural projection 
\[\bstar p:\gcov{\mathcal{L}}\;\longrightarrow \;\mathcal{L}\quad\quad\quad\quad
\bstar \tilde{x}\;\longmapsto\; \hat{x}\;=\;\lim p(\tilde{x}_{\alpha}) ,\]
where $\{ \tilde{x}_{\alpha}\}$ is a representative sequence
in the class $\bstar \tilde{x}$.  We will write $\lim \bstar\tilde{x}=\hat{x}$
if $\bstar p(\bstar\tilde{x})=\hat{x}$.  Note that $\bstar p$ is surjective
if and only if $L$ is dense and in general $\bstar p$ maps onto
the closure $\overline{L}$ of $L$, itself a sublamination of $\mathcal{L}$.

\begin{prop}\label{unionofgalaxies} $\gcov{\mathcal{L}}$ consists of a union of galaxies of 
$\;\bstar \widetilde{L}$.
\end{prop}

\begin{proof}  Let $\bstar\tilde{x}\in \gcov{\mathcal{L}}$ and denote
by $U$ the galaxy containing $\bstar\tilde{x}$.
If $\bstar\tilde{y}\in U$, then there exists a sequence of geodesic 
paths $\{\tilde{\eta}_{\alpha}\}$
connecting representatives $\{ \tilde{x}_{\alpha}\}$ to $\{ \tilde{y}_{\alpha}\}$
in $\widetilde{L}$, whose projection to $\mathcal{L}$ gives 
a convergent sequence of paths $\{\eta_{\alpha}\}$.
It follows that the projection $\{ p(\tilde{y}_{\alpha})\}$ converges, and
$\bstar\tilde{y}\in\gcov{\mathcal{L}}$ as well.
\end{proof}

The galaxies that make up $\gcov{\mathcal{L}}$ will be referred to as leaves.
See \cite{Ge1} for a proof of the following

\begin{theo}  $\gcom{\mathcal{L}}$ may be given the structure of a lamination
whose leaves are nowhere dense and such that the map 
$\bstar p: \gcom{\mathcal{L}}\rightarrow \overline{L}$ is an open surjection.
\end{theo}

One can thus think of $\gcov{\mathcal{L}}$ as a 
the result of unwrapping all of the diophantine approximations implied
by $L$.  The topology that $\gcov{\mathcal{L}}$ obtains from its lamination
atlas is not unique,  and is called a {\em germ universal cover
topology}.  It is in general coarser than the topology $\gcom{\mathcal{L}}$
induces from $\bstar \widetilde{L}$.

\begin{prop}\label{compactness}  If $\mathcal{L}$ is compact then 
$\gcov{\mathcal{L}}=\bstar\widetilde{L}$.
\end{prop}

\begin{proof}  This follows from well-known compactness arguments
{\em e.g.} see the proof in \cite{Ge1}.
\end{proof}

An element 
$\bstar u=\bstar f\circ\bstar g^{-1}\in \fgrm (\mathcal{L} ,L, x, {\sf f})$ 
arises as the limit of a sequence of $K_{\alpha}$-quasiisometries
\begin{equation}\label{quasiisometry}
 \Big\{ u_{\alpha} :D_{r_{\alpha}}(\tilde{x}_{\alpha})
\;\longrightarrow\; D_{r_{\alpha}}(\tilde{y}_{\alpha})\Big\} ,
\end{equation}
where $\{ \tilde{x}_{\alpha}\} ,\{ \tilde{y}_{\alpha}\}\subset L$, $K_{\alpha}\rightarrow 1$ 
and $r_{\alpha}\rightarrow\infty$.  The limit
$\bstar u :V\rightarrow U$
is independent of the sequence $\{ u_{\alpha}\}$
and depends only on the sequences of frames $\{{\sf f}_{x_{\alpha}}\}$,
$\{{\sf f}_{y_{\alpha}}\}$.
In particular, we could have obtained $\bstar u$ through 
the same sequence of quasiisometries with domains extended to 
a sequence of larger disks
$D_{s_{\alpha}}(\tilde{x}_{\alpha})$,
$s_{\alpha}>r_{\alpha}$ -- provided that the new quasiisometry constants 
converge to 1 as well.

Now for arbitrary $\bstar \tilde{x}\in\gcov{\mathcal{L}}$, 
the expression 
$\bstar u(\bstar x )$ does not even make formal sense, since 
$\bstar u$ is so far only defined
on the galaxy $V$.  
We contrast this with the constant curvature case, where,
because $\fgrm (\mathcal{L} , x,L, {\sf f})\subset {\sf Isom}(\bstar \M )$,
$\bstar u(\bstar\tilde{x})$ is always formally defined, although it need
not define
an element of $\gcov{\mathcal{L}}$.

Let us say that $\bstar u$ is {\em formally
defined} on an element $\bstar\tilde{w}\in\gcov{\mathcal{L}}$ if 
there exists a sequence 
(\ref{quasiisometry}) giving rise to $\bstar u$ and
a representative sequence $\{ \tilde{w}_{\alpha}\}$ of $\bstar\tilde{w}$ 
such that 
\[ D_{r'_{\alpha}}(w_{\alpha})\subset D_{r_{\alpha}}(\tilde{x}_{\alpha})\]
for all $\alpha$,
where $r'_{\alpha}\rightarrow\infty$.  It follows then that if $V'$ is
the galaxy containing $\bstar\tilde{w}$, then the limit $\bstar u$
is defined on $V'$ as well.  Whenever we write $\bstar u(\bstar w )$,
it will tacitly be understood that $\bstar u$ is formally defined
at $\bstar w$.

Define the domain of $\bstar u$ as
\[ {\sf Dom}(\bstar u)\;\; =\;\; 
\bigg\{ \bstar\tilde{x}\;\in\;\gcov{\mathcal{L}}\; \bigg|\;\;
\bstar u(\bstar\tilde{x})
\;\in\;\gcov{\mathcal{L}}\;\text{ and }\;
\lim \bstar u(\bstar\tilde{x} )\;=\;
\lim \bstar\tilde{x} 
\bigg\} ,\]
and ${\sf Ran}(\bstar u)=\bstar u({\sf Dom}(\bstar u))$.
With this definition, it follows that $\fgrm (\mathcal{L}, L, x, {\sf f})$
has the structure of a groupoid.  
Note that for any $\bstar u\in \fgrm (\mathcal{L},L, x, {\sf f})$,
${\sf Dom}(\bstar u )$, ${\sf Ran}(\bstar u)$ are unions of leaves and
hence induce lamination structures from
$\gcov{\mathcal{L}}$.  Moreover,
on ${\sf Dom}(\bstar u )$,
\begin{equation}\label{deckness} \bstar p\circ\bstar u\;\;=\;\;\bstar p .  
\end{equation}
In particular we see that 
$\bstar u : {\sf Dom}(\bstar u )\rightarrow {\sf Ran}(\bstar u )$
defines a lamination homeomorphism.

\begin{exam}  Let $\mathcal{L}$ be the irrational foliation of $\T^{2}$
by lines of slope $r$, $L\approx \R$ any dense leaf.  Then by 
Proposition~\ref{compactness}, $\gcov{\mathcal{L}}\;=\; \bstar\R$.
Moreover, for any frame field ${\sf f}$ and $\bstar u\in
\fgrm (\mathcal{L}, L,x, {\sf f})$, it is not difficult to see that
${\sf Dom}(\bstar u )=\bstar\R$.  Thus, $\fgrm (\mathcal{L}, L,x, {\sf f})$
is a group isomorphic to $\bast\Z_{r}$.
\end{exam}

\begin{exam} Let $\mathcal{L}$ be the profinite hyperbolic surface solenoid
$\widehat{\Sigma}_{\mathcal{G}}$ of {\em Example}~\ref{invlimitsol}.  Then we have, again
by compactness,
\[  \gcov{\mathcal{L}}\;\; =\;\; \bstar\HP^{2} . \]
If ${\sf f}$ is a frame field lifted from a frame on a surface occurring
in the defining inverse limit, then 
$\fgrm (\mathcal{L} ,L, x, {\sf f})$ is a group.  On the other hand, if ${\sf f}$ is a frame field
not obtained in this way, then 
$\fgrm (\mathcal{L} ,L, x, {\sf f})$ need not be a group {\em e.g.}\ see
\S 6.
\end{exam}

The following may also be found in \cite{Ge1}.

\begin{theo}\label{maplifting}  Let 
$F:(\mathcal{L},L)\rightarrow (\mathcal{L}',L')$ be a lamination
map.  Then there exist germ universal cover topologies so that
the map 
\[ \gcom{F}: \gcov{\mathcal{L}}\longrightarrow \gcov{\mathcal{L}'}\]
induced by $\{ \tilde{x}_{\alpha}\}\mapsto \{ \widetilde{F}(\tilde{x}_{\alpha})\}$
is a continuous lamination map.
\end{theo}

\begin{note}  It is useful here to point out that for a lamination 
$\mathcal{L}_{\rho} = ( \widetilde{B}\times {\sf F}) / \pi_{1}B$ 
occurring as a suspension
of a representation $\rho :\pi_{1}B\rightarrow {\sf Homeo}({\sf F})$,
it is in general false that a lamination homeomorphism
$F:\mathcal{L}_{\rho}\rightarrow \mathcal{L}_{\rho}$ lifts to a homeomorphism of
the ``universal covering space'' $\widetilde{B}\times {\sf F}$.
\end{note}

Now suppose that $L'$ is another  leaf of $\mathcal{L}$.  Denote by
$\gcov{\mathcal{L}}'$ the germ universal cover formed from $L'$.

\begin{prop}\label{virisomgrmcovers}  If $L'$ accumulates on $L$ then
there is a virtual subisometry 
$\bstar \widetilde{L}\rightarrow\bstar \widetilde{L}'$ restricting
 to a virtual subisometry
\[   \gcov{\mathcal{L}}\;\;  \longrightarrow \;\;\gcov{\mathcal{L}}'    \]
which is a homeomorphism onto its image with respect to appropriate
germ universal cover topologies.
\end{prop}

\begin{proof}  This follows directly from the proof of 
Theorem~\ref{leavesvirtrelated}.  
\end{proof}

\section{Sensitivity to Changes in Data}

In this section we shall examine the dependence of the fundamental germ
on the base point $x$, the accumulating leaf $L$ and the frame field ${\sf f}$.

\vspace{3mm}
\noindent
{\bf Change in base point and accumulating leaf:}  Let us fix for the moment the dense leaf
$L$ and consider a change of base point $x\mapsto x'$ in which $L_{x}=L_{x'}$.  
Let $\eta$ be a
geodesic connecting $x$ to $x'$ in $L_{x}$.  The tangent vector ${\sf v}$
to $\eta$ at $x$
has coordinate $(a_{1},\dots ,a_{n})$ with respect to the frame
${\sf f}_{x}$.  At each $y\in T =$ the domain of ${\sf f}$, 
this coordinate determines a vector ${\sf v}_{y}$
using the frame ${\sf f}_{y}$.  We obtain in this way
a transversally continuous family of geodesics $\{ \eta_{y}\}_{y\in T}$.
Restricting to an open subtransversal of $T$ if necessary,
we may parallel translate ${\sf f}$ along the geodesic family
to obtain a frame field ${\sf f}'$ with domain $T'\ni x'$.
The following is then immediate from the definition of the fundamental germ.

\begin{prop}  Let $x'$, ${\sf f}'$ be as in the preceding paragraph.
Then
\[  \fgrm (\mathcal{L},L,x,{\sf f} )\;\; =\;\;\fgrm (\mathcal{L},L,x',{\sf f}' ) .  \]
\end{prop}

If we consider a change of base point
$x\mapsto x'$, in which $L_{x}\not= L_{x'}$, 
the situation becomes considerably more 
subtle.  In fact, we shall see in \S\ref{antenna} that
fundamental germs based at points on different leaves can be
nonisomorphic.  For similar reasons, a change in accumulating leaf $L$
may yield nonisomorphic fundamental germs.

\vspace{3mm}
\noindent
{\bf Change in frame field:}  Let us now fix the base point $x$ and
consider a new frame field ${\sf f}': T'\rightarrow {\bf F}_{\ast}\mathcal{L}$
based at $x$.  For simplicity, we again assume that $\pi_{1}L=1$.
Since $T$ (the domain of ${\sf f}$) and $T'$ each contain subtransversal
neighborhoods
of $x$ lying in a common flow box, 
it is clear that there is a natural bijection
\[ \fgrm (\mathcal{L},L,x, {\sf f}) \;\longleftrightarrow \;
 \fgrm (\mathcal{L},L,x, {\sf f}'). \]
The issue is then the law of composition.  We will show that
this map need not be an isomorphism.

Let us consider the inverse limit solenoid $\widehat{\Sigma}_{\mathcal{G}}$ 
of {\em Example}~\ref{invlimitsol}. Assume that $L=L_{x}$, 
$T=T'=$ a fiber over a point $x_{0}\in\Sigma_{0}$ and that 
${\sf f}_{x}={\sf f}'_{x}$.  We will take
${\sf f}$ to be simply the lift of a frame based at $x_{0}$,
so that ${\sf f}$-diophantine approximations consist of sequences 
$\{\gamma_{\alpha}\}\subset\Gamma_{0}$ converging with respect
to the family $\{ \Gamma_{i}\}$.  It follows that every
${\sf f}'$-diophantine approximation of $x$ may be written in the form
$\big\{ \gamma_{\alpha}\Theta_{\alpha}  \big\}  $,
where $\{ \Theta_{\alpha}\}$ consists of a sequence of rotations based at $x$ with
angle going to 0 and $\{\gamma_{\alpha}\}$ is an ${\sf f}$-diophantine approximation.  
General elements of $\fgrm (\mathcal{L},L,x,{\sf f}')$
are then of the form
$  \big\{ \gamma_{\alpha}\Theta_{\alpha}\eta_{\alpha}  \big\} $,
where $\{\eta_{\alpha}\}$ is another ${\sf f}$-diophantine approximation.
We should not expect products of elements of this 
type to yield elements of $\fgrm (\mathcal{L},L,x,{\sf f}')$.  Indeed,
such a product would have the shape
\begin{equation}\label{product}  \big\{ \gamma_{\alpha}\Theta_{\alpha}
\eta_{\alpha}\Delta_{\alpha}\omega_{\alpha}  \big\}, \end{equation}
for $\{\Delta_{\alpha}\}$ another sequence of rotations based at $x$ with
angle going to $0$ and $\{\omega_{\alpha}\}$ another 
${\sf f}$-diophantine approximation.
If $\Theta_{\alpha}$ does not converge to the identity fast enough,
$\Theta_{\alpha}
\eta_{\alpha}\Delta_{\alpha}\omega_{\alpha}$ applied to ${\sf f}'_{x}$
will not project to a frame based at $x_{0}\in \Sigma_{0}$.
Hence the expression (\ref{product}) is not even asymptotic
to an element of $\fgrm (\mathcal{L},L,x,{\sf f}')$.
It is not difficult to see that unless ${\sf f}'$
is the pull-back of a frame on $\Sigma_{0}$,
this sort of problem always arises.

\section{The Mother Germ}\label{mother}

In this section, we assume that $\mathcal{L}$ has a dense
leaf $L$, with which we define the germ universal cover $\gcov{\mathcal{L}}$,
equipped with a fixed germ universal cover topology.

While the fundamental germ $\fgrm (\mathcal{L},L,x,{\sf f} )$ enjoys
the property of being reasonably calculable and leaf
specific, it can be sensitive
to data variation. There are additional shortcomings:
\begin{itemize}
\item By (\ref{deckness}), the action of the fundamental germ 
$\fgrm (\mathcal{L},L,x, {\sf f})$ on $\gcov{\mathcal{L}}$ 
respects the germ covering $\bstar p$. 
However it need not be the case that every identification implied by
$\bstar p$ is implemented by an element of $\fgrm (\mathcal{L},L,x, {\sf f})$.
\item There will be in general other 
maps of leaves of $\gcov{\mathcal{L}}$
that satisfy (\ref{deckness}) but do not appear in $\fgrm (\mathcal{L},L,x, {\sf f})$.
\item It appears that $\fgrm (\mathcal{L},L,x, {\sf f})$ such
as it is defined, will be functorial only under certain types
of lamination maps {\em e.g.}\ see \cite{Ge1}.
\end{itemize}
For this reason, we will expand 
$\fgrm (\mathcal{L},L,x, {\sf f})$ to a larger groupoid, called the mother germ. 
The mother germ will be the maximal amplification of 
$\fgrm (\mathcal{L},L,x, {\sf f})$ which contains all partially defined maps of
sublaminations of
$\gcov{\mathcal{L}}$ satisfying (\ref{deckness}): in other words, it is  
the full deck groupoid of $\bstar p$.

Let ${\sf Dom},\, {\sf Ran}$ be sublaminations of $\gcov{\mathcal{L}}$.
A homeomorphism
\[\bstar u :{\sf Dom}\;\longrightarrow \;{\sf Ran}\;\;\subset\;\;\gcov{\mathcal{L}}\]
satisfying (\ref{deckness}) is called {\em deck}.
Note that condition (\ref{deckness}) implies that a 
deck homeomorphism $\bstar u$ is automatically an isometry along
the leaves of ${\sf Dom}$.

\begin{defi}  The {\bf mother germ} is the groupoid
\[ \fgrm (\mathcal{L})\;\; =\;\; 
\Big\{\bstar u:{\sf Dom}\;\longrightarrow \;{\sf Ran} 
\text{ is a deck homeomorphism }\Big\}.   \]
\end{defi}

The mother germ will never be a group, since it distinguishes deck
maps obtained from others by restriction of domain.  In general, however,
it will contain many interesting and calculable subgroups and
subgroupoids, as
the following shows.

\begin{prop}\label{allinclude}  Let $L'$ be any dense leaf of $\mathcal{L}$. 
 Then there is an injective groupoid homomorphism
\[  \fgrm (\mathcal{L},L',x, {\sf f})\;\hookrightarrow\; 
\fgrm (\mathcal{L}) .  \]
\end{prop}

\begin{proof}  By Proposition~\ref{virisomgrmcovers}, there exists
an isometric inclusion 
$\bstar f:\gcov{\mathcal{L}}'\hookrightarrow\gcov{\mathcal{L}}$.  
If $\bstar u\in\fgrm (\mathcal{L},L',x, {\sf f})$,
then the map 
\[ \bstar u\;\longmapsto\; \bstar f\circ \bstar u\circ\bstar f^{-1}\] 
defines an injective groupoid homomorphism. 
\end{proof}

\begin{theo}\label{quotient}  The quotient 
\[  \fgrm (\mathcal{L})\big\backslash \gcov{\mathcal{L}} ,  \]
equipped with the quotient germ universal cover topology,
has the structure of a riemannian lamination canonically isometric to $\mathcal{L}$.
\end{theo}

\begin{proof}  Let $\bstar \tilde{x}$, $\bstar \tilde{y}\in\gcov{\mathcal{L}}$
be such that $\lim\bstar \tilde{x}  =
\lim  \bstar \tilde{y} =x$.  Thus each point is represented
by sequences in $\widetilde{L}$ that project to sequences $\{ x_{\alpha}\}$,
$\{ y_{\alpha}\}\subset L$ having a common limit $x$.  Let ${\sf f}$
be a frame field along a transversal $T$ containing $x$ and which we may
assume contains $\{ x_{\alpha}\}$ and
$\{ y_{\alpha}\}$.  Then if $\bstar f$, $\bstar g$ are the 
diophantine approximations associated to $\{ x_{\alpha}\}$,
$\{ y_{\alpha}\}$ we have $\bstar u = \bstar g\circ \bstar f^{-1}\in
\fgrm (\mathcal{L},L,x, {\sf f})$
identifies $\bstar \tilde{x}$ with $\bstar \tilde{y}$.  Since this latter
groupoid belongs to the mother germ by Proposition~\ref{allinclude}, it follows
that $ \fgrm (\mathcal{L})\backslash \gcov{\mathcal{L}}$ contains all
of the identifications implied by $\bstar p$ and so may be identified
with $\mathcal{L}$ with its quotient topology.
\end{proof}

Let $\mathcal{L}$, $\mathcal{L}'$ be riemannian laminations with dense
leaves $L, L'$.  A map 
$\gcom{F}:\gcov{\mathcal{L}}\rightarrow \gzg{\widetilde{\mathcal{L}}'}$
is called {\em standard} if it is induced by 
$F:\widetilde{L}\rightarrow \widetilde{L}'$ ({\em e.g.}\ compare with the 
definition found in \S 3).
In addition, $\gcom{F}$ is called {\em $\fgrm (\mathcal{L})$-equivariant} 
if there exists
a groupoid homomorphism $\gzg{F}_{\ast}:\fgrm (\mathcal{L})\rightarrow
\fgrm (\mathcal{L}')$ such that
\[ \gcom{F}\big(\bstar u\cdot \bstar \tilde{x} \big) 
\;\; =\;\; 
\gzg{F}_{\ast}\big(\bstar u\big)\cdot \gcom{F}\big(\bstar \tilde{x} \big)  \]
for all $\bstar u\in \fgrm (\mathcal{L})$ and 
$\bstar \tilde{x}\in {\sf Dom}(\bstar u)$.

\begin{theo}\label{descend}  Let 
$\gcom{F}:\gcov{\mathcal{L}}\rightarrow \gzg{\widetilde{\mathcal{L}}'}$ 
be a standard, $\fgrm (\mathcal{L})$-equivariant map.  Then
$\gcom{F}$ covers a unique map $F:\mathcal{L}\rightarrow\mathcal{L}'$.
\end{theo}

\begin{proof}  By equivariance, the expression 
\[   F  \;\;=\;\; \bstar p'\circ \gcom{F}\circ\bstar p^{-1}  \]
yields a well-defined function
$F:\mathcal{L}\rightarrow\mathcal{L}'$, continuous because
$\bstar p ,\bstar p'$ are open maps and $\gcom{F}$ is continuous. 
\end{proof}

In \cite{Ge1}, functoriality of the fundamental germ was 
demonstrated only with respect to the restricted class of {\em trained} 
lamination maps.  The following theorem shows that
the mother germ is considerably more flexable.
A {\em lamination covering map} is a surjective lamination map
which is a covering map when restricted to any leaf.

\begin{theo}\label{functor}  Let $F: \mathcal{L}\rightarrow \mathcal{L}'$
be a lamination covering map.  Then $F$ induces
an injective homomorphism of mother germs
\[ \gzg{F}_{\ast}:\fgrm (\mathcal{L})\;\hookrightarrow\; 
\fgrm (\mathcal{L}'). \]
\end{theo}

\begin{proof}  Let $L\subset\mathcal{L}$ be a dense leaf and let  
$\widetilde{F}:\widetilde{L}\rightarrow\widetilde{L}'$
be the leaf universal cover lift.  Then by Theorem~\ref{maplifting},
$\widetilde{F}$ induces
a standard map
\[\gcom{F}:\gcov{\mathcal{L}}\longrightarrow\gcov{\mathcal{L}'}.\]  
We note that since $\widetilde{F}$
is injective, $\gcom{F}$ is a homeomorphism onto its image lamination.
Let
$\bstar u\in\fgrm (\mathcal{L})$.  Then the map
\[  \gzg{F}_{\ast}\big(\bstar u\big)\;\; :=\;\;
\gcom{F}\circ\bstar u\circ \gcom{F}^{-1} ,   \]
defined on $\gcom{F}\big({\sf Dom}(\bstar u)\big)$,
is deck for the germ universal covering $\bstar p'$.  Indeed
\begin{eqnarray*}  \bstar p'\circ
\Big(  \gcom{F}\circ\bstar u\circ  
\gcom{F}^{-1} \Big)  & =& F\circ 
(\bstar p\circ\bstar u) \circ \gcom{F}^{-1}\\
& =& F\circ \bstar p\circ \gcom{F}^{-1}\\
& =& \bstar p'.
\end{eqnarray*}
Thus the map $\gzg{F}_{\ast}$ is an injective groupoid homomorphism, and we are done.
\end{proof}
 
We have the following 

\begin{coro}  The mother germ $\fgrm (\mathcal{L})$ is independent of
leaf-wise riemannian metric and smooth structure.  In particular, 
$\fgrm (\mathcal{L})$
is a {\rm topological} invariant.
\end{coro}

\section{The Antenna Lamination}\label{antenna}

In this section, we will calculate the fundamental germ of the
antenna Riemann surface lamination of Kenyon and Ghys \cite{Gh}:
it is distinguished by the unusual property of having dense leaves of both
planar and hyperbolic conformal type.  

We begin by constructing a graphical model of a dense leaf 
of the antenna lamination.  
Let $T_{1}$ be the cross with vertices $V_{1}=\big\{
(0,0), (\pm 1 ,0), (0,\pm 2) \big\}$ and edges consisting of the
line segments connecting $(0,0)$ to each of the other four
vertices.  Suppose that we have constructed $T_{n}$ meeting
the $x$-axis in the interval $[-2^{n}+1,2^{n}-1]\times\{ 0\}$ and
meeting the $y$-axis in the interval $\{ 0\}\times
[-2^{n},2^{n}]$.  Translate $T_{n}$ vertically so that the origin
is taken to $(0,2^{n})$ and consider the images of this translate
by rotations of the plane -- about the origin -- of angles $0,\pm
\pi/2,\pi$.  The union of these images forms a tree; $T_{n+1}$ is
then obtained by replacing the extremal edges 
$[2^{n+1}-2,2^{n+1}]\times\{ 0\}$ and $[-2^{n+1},-2^{n+1}+2]\times\{ 0\}$
by  $[2^{n+1}-2,
2^{n+1}-1]\times\{ 0\}$ and $[-2^{n+1}+1,-2^{n+1}+2]\times\{ 0\}$.  It follows that $T_{1}\subset
T_{2}\subset\dots$: we then define
\[ T_{\infty}\;\; =\;\;\lim_{\longrightarrow}T_{n}.\]
See Figure 1.

\begin{figure}[htbp]
\centering \epsfig{file=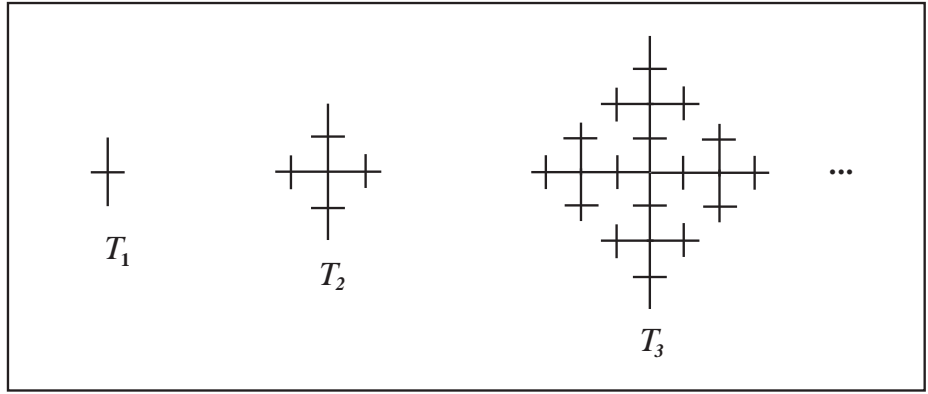, width=5in} \caption{The
Antenna Tree}
\end{figure}

Given $n\in \Z$, let ${\sf ord}_{2}(n)$ be the 2-adic order: 
the largest nonnegative integer $r$ for which $2^{r}$ divides $n$.
Then the vertex
set of $T_{\infty}$ is
\[ V_{\infty} \;\; =\;\; \big\{ (0,0)\big\}\;\bigcup\;
\Big\{ v=(x,y)\in\Z\oplus\Z\;\big|\;\; {\sf ord}_{2}(x)\not= {\sf ord}_{2}(y) \Big\}.\]
We may view $V_{\infty}$ as a groupoid through its action on itself 
by addition.  In order to avoid confusion, 
we write $v\circ w$ to indicate groupoid
composition, in order to distinguish it from
the element $v+w\in\Z\oplus\Z$. 

\begin{prop}\label{vertexgpd}  For all $v,w\in V_{\infty}$,
the composition $v\circ w$ is defined if and only if $v=-w$.
\end{prop}

\begin{proof}   Let $v,w\in V_{\infty}$.  We show that
${\sf Ran}(w)={\sf Dom}(v)$ if and only if $v=-w$. Suppose $v\not= -w$.
Then we may write 
\[ v+w\;\;=\;\; \left( \sum_{\alpha=m}^{M}a_{\alpha}2^{\alpha},
\sum_{\alpha=n}^{N}b_{\alpha}2^{\alpha}\right)  \]
where $m$, $n$ are the first non-zero indices of the $2$-adic expansions
of the coordinates.  If $v\circ w$ is defined, then since $0\in {\sf Dom}(w)$,
we must have $v+w\in V_{\infty}$.  In particular,
at least one of $m$ or $n$ is nonzero.  Suppose it is $n$; we may
assume without loss of generality that $m<n$.  Write
\[  w\;\;=\;\; \left( \sum_{\alpha=r}^{R}c_{\alpha}2^{\alpha},
\sum_{\alpha=s}^{S}d_{\alpha}2^{\alpha}\right) .\]  
Let $x=(0,2^{m})$.  If $r<m$,
then $x\in {\sf Dom}(w)$ but $v+w+x\not\in V_{\infty}$ {\em i.e.}
${\sf Ran}(w)\not= {\sf Dom}(v)$.  This is also true
when $r\geq m$ except for two cases.
If $r>m$ and $s=m$, 
$w+x$ is not defined presisely
when $1=d_{s}=\dots = d_{r-1}$ and $d_{r}=0$. Here we take 
$x'=(2^{r},2^{m})\in {\sf Dom}(w)$ and note that
$v+w+x'\not\in V_{\infty}$. If $r=m$ and $s>m$, then $w+x$
is not defined.  In this case, it follows from the form of $v+w$
that if $v=(v_{1},v_{2})$ then ${\sf ord}_{2}(v_{1})>m$, so that $x\in {\sf Dom}(v)$.
On the other hand, $x-w\notin {\sf Dom}(w)$.  Thus
${\sf Ran}(w)\not={\sf Dom}(v)$ here as well.  
\end{proof}

The lines $x=\pm y$ intersect $T_{\infty}$ at the
origin only.  Each of the four components of $T_{\infty}\setminus
\{(0,0)\}$ defines an end, one contained in each of the four components
of $\R^{2}\setminus \{ (x,\pm x)\; |\;\;x\in\R \}$. Equipped with 
the path metric induced from $\R^{2}$, 
$T_{\infty}$ has exactly four orientation preserving isometries, corresponding to
the rotations about the origin of angles $0,\pm
\pi/2,\pi$ (since ends must be taken to ends).  On the other hand,
$T_{\infty}$ has many partially defined isometries.  For example, for 
$v\in V_{\infty}$, let $I_{v}$ be the map of $\Z\oplus\Z$ defined 
\[ I_{v}(x,y)\;\; =\;\; v+(x,y) .\] Then there is a maximal
subtree $T_{\infty}^{v}\subset T_{\infty}$ (not necessarily connected) for
which $I_{v}(T_{\infty}^{v})\subset T_{\infty}$.  By definition, $I_{v}$ is
isometric on its domain of definition.  
If $v$ has coordinates of large 2-adic order, then $I_{v}$
is defined on a large ball about $0$ in $T_{\infty}$.  More precisely,
if $v=(x,y)$ and ${\sf ord}_{2}(x),{\sf ord}_{2}(y)\geq n$ then
$T_{n}\subset T_{\infty}^{v}$.
Although
the inverse $I_{v}^{-1}=I_{-v}$ is always defined at $0$,
the composition $I_{v_{1}}\circ I_{v_{2}}=I_{v_{1}+ v_{2}}$ will not be defined 
at $0$ if
$v_{1}+v_{2}\notin V_{\infty}$.

We now define a riemannian surface modelled on $T_{\infty}$, which
will occur as a dense leaf of the antenna lamination.
Regarding $T_{\infty}\subset\R^{2}\times \{ 0\}\subset\R^{3}$, it is clear
that
\[ S_{\infty}\;\; = \;\; \text{boundary of a tubular neighborhood of }T_{\infty}\] 
is homeomorphic to a sphere with four punctures.  We want to fix a particular 
realization of $S_{\infty}$ so that the partial isometries 
$I_{v}$ of $T_{\infty}$ will induce partial isometries of $S_{\infty}$. 
Torward this end, consider the surfaces shown in Figure 2. 
\begin{figure}[htbp]
\centering \epsfig{file=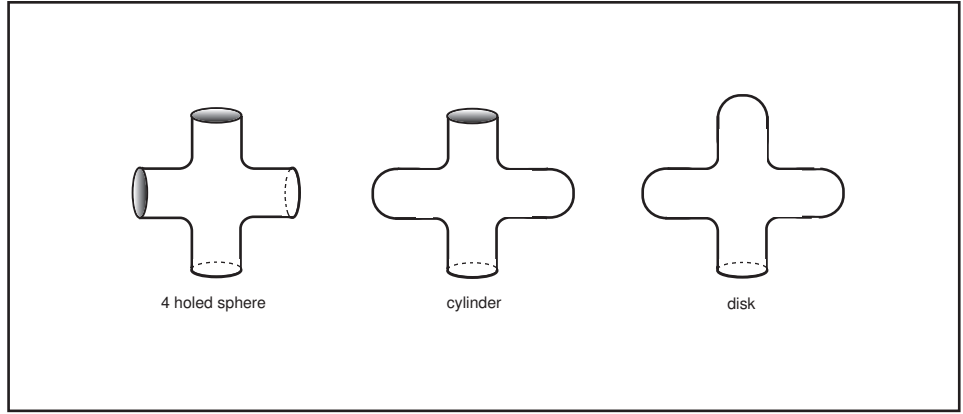, width=5in} \caption{Building Blocks}
\end{figure}
We assume that they are equipped with riemannian metrics and
boundary parametrizations so that given any pair of such surfaces and
a choice of boundary component of each, the glueings are canonical
and isometric. Each riemannian surface corresponds to a subgraph of $T_{\infty}$,
and we may build $S_{\infty}$ from these riemannian surfaces using $T_{\infty}$ 
as a template.  The metrics on the building blocks will also 
be chosen
so that when $S_{\infty}$ is assembled within $\R^{3}$ it is invariant
not only with respect to $\pi /2$-rotations about the $z$-axis, but also
$\pi$-rotations about the $x$ and $y$-axes. 
We think of $T_{\infty}$ as a spine floating inside 
the tubular neighborhood bounded by $S_{\infty}$, and 
we project in the positive vertical direction a copy of $T_{\infty}$ onto $S_{\infty}$. 
We denote this copy also by $T_{\infty}$, and use the symbols $0$ and $v$ to
denote the origin and a generic element of its vertex set as well.
Having constructed $S_{\infty}$ in this way, 
it is clear that every $I_{v}$ induces
a partial isometry of $S_{\infty}$ whose domain is the subsurface (with boundary) of $S_{\infty}$ modelled
on $T_{\infty}^{v}$.  We denote this partial isometry $I_{v}$ as well.

Let $S_{\infty}^{+}$ be the intersection of $S_{\infty}$ with the half plane
$z\geq 0$.  The universal
cover $\widetilde{S}_{\infty}$ of $S_{\infty}$ is built up from ``tiles''  modelled
on $S_{\infty}^{+}$, glued together side by side according to the same pattern
one uses to glue ideal quadrilaterals to obtain the hyperbolic plane as 
the universal cover of the four times punctured sphere.
Fix 
$\tilde{0}\in \widetilde{S}_{\infty}$ a base point
lying over $0$.  The deck group of the universal covering map is $F_{3}$, the
free group on three generators.

Let ${\sf w}_{0}$ be the unit vector based $0$ which is
parallel to the $x$-axis and points in the positive direction. 
Consider the vector field ${\sf W}$ on the vertices of $T_{\infty}$ 
obtained by parallel translating ${\sf w}_{0}$ along $T_{\infty}$.  
Note that the partial isometry of $S_{\infty}$ induced by $I_{v}$, 
$v\in V_{\infty}$, takes ${\sf w}_{0}$ to ${\sf w}_{v}={\sf W}(v)$.  
This is not true of the rotations by angles $\pm \pi /2$ and $\pi$.

Let ${\sf D}_{\tilde{0}}\subset \widetilde{S}_{\infty}$ be the fundamental domain 
containing $\tilde{0}$.  We lift $T_{\infty}$ to ${\sf D}_{\tilde{0}}$, then translate
it by $F_{3}$ to obtain a (disconnected) graph 
$\widetilde{T}_{\infty}$ on $\widetilde{S}_{\infty}$.
Let $\widetilde{\sf W}$ be the vector field defined on the vertices of $\widetilde{T}_{\infty}$
that is the lift of ${\sf W}$.
The partial isometry $I_{v}$ lifts to a partial isometry of $\widetilde{S}_{\infty}$
which maps a region of each fundamental domain ${\sf D}$ into ${\sf D}$: we denote this privileged lift by $I_{v}$ 
as well, and the set of such privileged lifts is denoted $\I$.
In addition, by composing with
elements of $F_{3}$, we obtain new partial isometries covering $I_{v}:S_{\infty}\rightarrow S_{\infty}$.
We denote by $\widetilde{\I}$ the set of partial isometries of $\widetilde{S}_{\infty}$
obtained in this way.    
Then $F_{3},\, \I\subset \widetilde{\I}$, and every element of $F_{3}$ commutes
with every element of $\I$.

We are now ready to describe the antenna lamination.  
Consider first the space $\A$ of all trees in $\R^{2}$ whose
vertex set contains the origin $0$ and lies within $\Z\oplus\Z$.  Each
tree $T\in\A$ is equipped with the path metric induced from $\R^{2}$.  On $\A$, we
consider the metric
\[ d(T,T')\;\; =\;\; \exp (-n), \]
where $n$ is the largest integer such that the ball of radius $n$ about
$0$ in $T$ coincides with that about $0$ in $T'$.  $\A$ is a
compact metric space, \cite{Gh}.  Two graphs $T$ and $T'$ are termed equivalent
if there exists a translation by
$(x,y)\in\Z\oplus\Z$ such that $T+(x,y)=T'$.

Now for any tree $T\in\A$, the ball of radius 1 about $0$ is a tree
$P\in\A$
all of whose vertices lie in the set $\{(0,0)\}\cup\{(\pm
1,\pm 1)\}$.  We write $|P|\leq 4$ for the number of vertices $v$ of $P$
different
from $0$.  There are 16 possible such $P$, and we may decompose $\A$
into a disjoint union of clopens $\A_{P}$, where $\A_{P}$  consists of
those trees whose unit ball
about $0$ is $P$.

For each $P$, we consider in the spirit of Figure 2
a model pointed Riemann surface $(\Sigma_{P},z_{P})$ 
homeomorphic to $\SI^{2}\setminus (|P|\text{ open
disks})$.  We assume as before that each boundary component
$\partial_{v}\Sigma_{P}$ -- labeled by a vertex $v\not= (0,0)$ of
$P$ -- has a fixed paramentrization, so that any two
may be identified along their boundaries without ambiguity. 
Define 
\[ \mathcal{L} \;\; =\;\;
\Big( \bigcup\, \big( \A_{P}\times\Sigma_{P}\big)\Big)\Big/\text{ gluing} ,
\] where the gluing is performed as follows.  Given $T\in\A_{P}$,
$v\in P$, the translate $T+v$ is in $\A_{P'}$ for some $P'$, where
$-v\in P'$.  We then glue the boundaries $\partial_{v}\Sigma_{P}$
and $\partial_{-v}\Sigma_{P'}$.  These gluings are compatible with
the trivial lamination structures on the $\A_{P}\times\Sigma_{P}$
and thus $\mathcal{L}$ has the structure of a riemannian surface
lamination. Note that there is an embedding
$\A\hookrightarrow\mathcal{L}$ induced by $\A_{P}\times\{
z_{P}\}\hookrightarrow\A_{P}\times\Sigma_{P}$.

Each leaf $L\subset \mathcal{L}$ corresponds to an equivalence
class of graph $T\in\A$, embedded in $L$ as a spine.  Note that
$S_{\infty}$ is the leaf corresponding to the class of
$T_{\infty}$. Define the {\em antenna lamination}
$\mathcal{L}_{\infty}$ to be the closure of
$S_{\infty}$ in $\mathcal{L}$.

Denote by $S_{n}\subset S_{\infty}$
the surface (with boundary) modelled on the subgraph $T_{n}\subset T_{\infty}$.  
If centered at a vertex $v\in V_{\infty}$ there is
a subgraph isometric to $T_{n}$, it models a subsurface $S_{n}(v)$ containing $v$,
and the isometry $I_{v}$ maps $S_{n}$ to $S_{n}(v)$.

The closure of
$V_{\infty}$ in $\mathcal{L}_{\infty}$ defines
a transversal ${\sf T}$ through $0\in S_{\infty}$, and the vector field 
${\sf W}$ is transversally
continuous with respect to the topology of ${\sf T}$. 
A point $v\in V_{\infty}$ is transversally close to $0$ if and
only if its coordinates have large 2-adic order.

We are now ready to calculate the fundamental germ 
\[\fgrm (\mathcal{L}_{\infty},0, {\sf f}),\] where ${\sf f}$
is the orthonormal frame field determined by ${\sf W}$. 

Define nested sets
\[  \widetilde{G}\; =\; \{ \widetilde{G}_{n}\}\;\;\subset\;\;
 \widetilde{\I}\;\;\;\;\text{ and }
\;\;\;\; G \; =\; \{ G_{n}\}\;\;\subset\;\; \I ,\]
$n=0,1,2,\dots $,
as follows.  We say that $\widetilde{I}\in\widetilde{\I}$ is $n$-{\it close} if the domain of
$\widetilde{I}$ contains the finite tree
$\widetilde{T}_{n}\subset\widetilde{T}_{\infty}\cap {\sf D}_{\tilde{0}}$
corresponding to $T_{n}$, and maps it into the
fundamental domain containing $\widetilde{I}(\tilde{0})$.  Then
$\widetilde{G}_{n}$
consists of the set
of $n$-close maps and $G_{n}$ the $n$-close maps in $\I$.
Observe that
\[ F_{3}\;\; = \;\; \bigcap \widetilde{G}_{n} . \]

For $I\in G_{n}$, $I^{-1}\in G_{n}$ also.  
Moreover, if $I'\in G_{m}$ and the composition $I\circ I'$
is defined at $\tilde{0}$, then it belongs to $G_{N}$,
for $N=\min
(m,n)$.

Let 
\[ \gzg{G} = \Big\{ \big\{ I_{v_{\alpha}}\circ I_{v'_{\alpha}}^{-1}\big\}\; =\;
\big\{ I_{v_{\alpha}-v'_{\alpha}}\big\}\;\Big|\;\; 
\big\{ I_{v_{\alpha}}\big\},\; \big\{ I_{v_{\alpha}'}\big\}\subset \I\text{ and converge w.r.t\ } G
\Big\}\Big/ \sim , \]
where the relation $\sim$ is defined by an ultrafilter $\mathfrak{U}$ as in
(\ref{nonstandardspace}).  We denote the elements of $\gzg{G}$ by $I_{\bast v}$
where $\bast v\in\bast\Z\oplus\bast \Z$, and regard $\gzg{G}$ as a groupoid with unit space
$\gcov{\mathcal{L}_{\infty}}$, in which the domains of elements are taken to be
maximal in the sense defined in \S 4.

\begin{prop}\label{complaw} As a set, $\gzg{G}$ may be identified with 
$\bast\Z_{\hat{2}}\oplus\bast\Z_{\hat{2}}$, where
\[ \bast\Z_{\hat{2}}\;\; =\;\; \Big\{ \{n_{\alpha}\}\subset\Z\;\Big|\;\;
{\sf ord}_{2}(n_{\alpha})\rightarrow \infty \text{ as }\alpha\rightarrow\infty  \Big\}\Big/\sim.  \]
Given $\bast v, \bast w\in\bast\Z_{\hat{2}}\oplus\bast\Z_{\hat{2}}$, 
the composition $I_{\bast v} \circ I_{\bast w}$ 
is defined if and only if $\bast v=-\bast w$.
\end{prop}

\begin{proof}  Any element 
$(\bast n_{1},\bast n_{2})\in \bast\Z_{\hat{2}}\oplus\bast\Z_{\hat{2}}$ 
may be written $(\bast n_{1},0)-(0,-\bast n_{2})$ which clearly defines
an element of $\gzg{G}$.  Now consider 
$\bast v, \bast w\in \bast\Z_{\hat{2}}\oplus\bast\Z_{\hat{2}}$, and
suppose that $\bast v\not=-\bast w$.  
The $2$-adic order extends to 
\[ {\sf ord}_{2}:\bast\Z_{\hat{2}}\rightarrow \bast\Z_{\infty}=(\bast\Z\setminus
\Z )\cup \{ 0\}. \]
Let $\bast V_{\infty}\subset \bast\Z_{\hat{2}}\oplus\bast\Z_{\hat{2}}$
be the subset of pairs $\bast u=(\bast u_{1},\bast u_{2})$ for which
${\sf ord}_{2}(\bast u_{1} )\not= {\sf ord}_{2}(\bast u_{2} )$.  We distinguish
four cases depending on whether $\bast v ,\bast w\in \bast V_{\infty}$ or not.
If $\bast v ,\bast w\in \bast V_{\infty}$ then we may regard each as a class
of sequences $\{ v_{\alpha}\}$, $\{ w_{\alpha}\}\subset V_{\infty}$.
If $\{ m_{\alpha}\}$, $\{ n_{\alpha}\}$ , $\{ r_{\alpha}\}$, 
$\{ s_{\alpha}\}$ are the 
sequences of indices occurring as in Proposition~\ref{vertexgpd}, then
there classes $\bast m$, $\bast n$, $\bast r$, $\bast s$ are totally
ordered in $\bast \Z$, hence we may assume the representative sequences
are.  In particular, we may proceed with the same argument as in
Proposition~\ref{vertexgpd}: the sequences $\{ x_{\alpha}\}$, $\{ x_{\alpha}'\}$
define elements of $\gcov{\mathcal{L}_{\infty}}$ which may be used to show that 
the composition $I_{\bast v} \circ I_{\bast w}$ is not defined.
Now suppose that $\bast v\notin\bast V_{\infty}$ but $\bast w\in\bast V_{\infty}$.
This means that both components of $\bast v$ have the same order denoted
${\sf ord}_{2}(\bast v )$.  
Then there exists $\bast x\in\bast V_{\infty}$ such that 
$\bast w+\bast x\in\bast V_{\infty}$
in which the two components of
$\bast w+\bast x$ have order greater than 
${\sf ord}_{2}(\bast v )$.  Then both components of $\bast v +\bast w+\bast x$
have equal order, which implies that $I_{\bast v} \circ I_{\bast w}$ is not defined.  
The case where
$\bast v\in\bast V_{\infty}$ but $\bast w\notin\bast V_{\infty}$ is handled
similarly.  Now suppose $\bast v,\bast w\notin\bast V_{\infty}$.  Here there
are two subcases.  First suppose that the orders of the components of
$\bast v$, $\bast w$ are not equal.  Denote by ${\sf ord}_{2}(\bast v )$,
${\sf ord}_{2}(\bast w )$ the common
order of the components of $\bast v$, $\bast w $.  Then
if say ${\sf ord}_{2}(\bast v ) <{\sf ord}_{2}(\bast w )$, 
we define $\bast x = (0,\bast w_{2})$
where $\bast w_{2}$ is the second component of $\bast w$.  Then $I_{\bast w}(\bast x)$
is defined but $I_{\bast v +\bast w}( \bast x)$ is not.  If
 ${\sf ord}_{2}(\bast v ) >{\sf ord}_{2}(\bast w )$ then $I_{\bast v}$ is defined
 on $\bast y= (0, \bast v_{2})$ but $\bast y-\bast w$ does not define
 an element of ${\sf Dom}(I_{\bast w} )$ since it does not converge
to the same point in $\gcov{\mathcal{L}_{\infty}}$ as $\bast y$. 
What remains is the case when ${\sf ord}_{2}(\bast v ) ={\sf ord}_{2}(\bast w )$.
If $\bast v +\bast w$ lies in $\bast V_{\infty}$ then
$\bast w\in {\sf Dom}(I_{\bast v} )$ but not in ${\sf Ran}(I_{\bast w} )$.
Otherwise, if the norms of the components of $\bast v +\bast w$
are equal, then ${\sf ord}_{2}(\bast v +\bast w)> {\sf ord}_{2}(\bast v )=
{\sf ord}_{2}(\bast w )$.
If we let $\bast x = ((\bast v+\bast w)_{1},0)$ then 
$\bast x\in {\sf Dom}(I_{\bast v+\bast w })$ but not in ${\sf Dom}(I_{\bast w} )$
so it cannot be that $I_{\bast v}\circ I_{\bast w}=I_{\bast v+\bast w}$.
\end{proof}

\begin{theo}\label{ellinfty}  As a set
\[ \fgrm (\mathcal{L}_{\infty},0,{\sf f})\;\; =\;\; 
\bast F_{3}\times \gzg{G}.\]  
The composition $\bstar u\circ \bstar v$, 
where $\bstar v =(\bast x,I_{\bast v})$ , $\bstar w=(\bast y ,I_{\bast w})$
is defined if and only if $\bast v=-\bast w$.
\end{theo}

\begin{proof} Every element $\widetilde{\I}$ may be written in
the form $I_{v}\circ\gamma =\gamma\circ I_{v}$ for $v\in V_{\infty}$
and $\gamma\in F_{3}$.  
Moreover, if $\widetilde{I}\in \widetilde{G}_{n}$, then $I\in G_{n}$.
The second statement follows immediately from Proposition~\ref{complaw}.
\end{proof}

Thus, although $\bast F_{3}\times \gzg{G}$ is formally a group,  
$\fgrm (\mathcal{L}_{\infty},0,{\sf f})$ is not a group with respect to 
the groupoid structure defined by its action on
$\gcov{\mathcal{L}_{\infty}}$.  It has nevertheless a distinguished
subgroup isomorphic to $\bast F_{3}$.
On the other hand, 

\begin{theo} Any two elements $\bstar v$
and $\bstar w$ of $\fgrm (\mathcal{L}_{\infty},0,{\sf f})$
define composable elements of the mother germ 
$\fgrm (\mathcal{L}_{\infty})$ by restriction of
domains. 
\end{theo}

\begin{proof}  Let $I_{\bast v}$, $I_{\bast w}$ be the $\gzg{G}$-coordinates 
of $\bstar v$, $\bstar w$.  If $\bast v$, $\bast w$ and 
$\bast v+\bast w$ belong to $\bast V_{\infty}$ then by restricting
$\bstar w$ to the leaf $S_{\infty}$ and restricting $\bstar v$
to the leaf of $\gcom{\mathcal{L}_{\infty}}$ containing $\bstar w$,
we obtain composable elements of $\fgrm (\mathcal{L}_{\infty})$. The other
cases are handled similarly and are left to the reader.
\end{proof}

The lamination $\mathcal{L}_{\infty}$ has the following
property: every leaf $L\not= S_{\infty}$ is
conformal to either $\C$ or $\C^{\ast}$ = $\C\setminus
\{(0,0)\}$, \cite{Gh}.  Hence $\mathcal{L}_{\infty}$ is neither a suspension nor a
locally free action of a Lie group.  In particular, the antenna lamination
is beyond the purview of the definition
of $\fgrm$ found in \cite{Ge1}.

Given any leaf $L$ of $\mathcal{L}_{\infty}$, one
can obtain a graphical model $T$ of $L$ as a limit of a sequence
of translations of $T_{\infty}$. One can then 
repeat the discussion leading up to Theorem~\ref{ellinfty} for $L$.
The proof of the following is left to the reader.

\begin{theo}  Let $L\subset\mathcal{L}_{\infty}$ be any leaf, modelled
as above on a graph $T\in\A$ with vertex set $V$.  Then for $v\in V$ and ${\sf f}$
constructed using a vector field as above, 
$\fgrm (\mathcal{L}_{\infty},v,{\sf f})$ may be identified with
\[ \bast\pi_{1}L\times \gzg{G},  \]
where $\bast\pi_{1}L\times \{ 0\}$ is a subgroup with respect to 
the groupoid structure that is $\cong$ $1$ or $\bast\Z$.
\end{theo}

\begin{coro}  Let $v\in V_{\infty}\subset S_{\infty}$ and $v'\in V'\subset L'\not=S_{\infty}$.  
Then choosing frame fields as above,
the fundamental germs $\fgrm (\mathcal{L}_{\infty},v,{\sf f})$ 
and $\fgrm (\mathcal{L}_{\infty},v',{\sf f}')$ are not isomorphic.
\end{coro}

\begin{proof}  $\fgrm (\mathcal{L}_{\infty},v,{\sf f})$ has a nonabelian
subgroup whereas $\fgrm (\mathcal{L}_{\infty},v',{\sf f}')$ is an abelian groupoid.
\end{proof}

\section{The $PSL(2,\Z )$ Anosov Foliation}

Let $\Gamma \subset PSL(2,\R )$ be a discrete group of finite type,
possibly with elliptic elements.  The quotient $\Sigma=\HP^{2}/\Gamma$ is a
finite volume hyperbolic surface orbifold.  The unit tangent bundle $T_{1}\Sigma$ is
defined to be the quotient ${\rm T}_{1}\HP^{2}/\Gamma$.
Let $\rho :\Gamma\rightarrow {\sf Homeo}(\SI^{1})$ be the representation
obtained by extending the action of $\Gamma$ to the boundary of $\HP^{2}$.
Then ${\rm T}_{1}\Sigma$ may be identified with
the suspension
\[ \Big( \HP^{2}\times \SI^{1}\Big)\Big/ \Gamma  \]
as follows.  Given  $(\tilde{z},t)\in \HP^{2}\times \SI^{1}$, associate
${\sf v}_{\tilde{z}}\in {\rm T}_{1}\HP^{2}$, the vector based at
$\tilde{z}$ and tangent to the ray limiting
to $t$.  This association is $\Gamma$-equivariant and descends to the
desired homeomorphism.  The expression of ${\rm T}_{1}\Sigma$ as a suspension defines
a hyperbolic Riemann surface foliation $\mathcal{F}$ on ${\rm T}_{1}\Sigma$,
which is also a fiber bundle over $\Sigma$ provided that $\Gamma$ has no
elliptic points. $\mathcal{F}$ is called an {\em Anosov foliation}.

In \cite{Ge1}, we worked with a definition of  $\fgrm$
that was available for suspensions such as $\mathcal{F}$ 
formed from fixed point free $\Gamma$.  Unfortunately,
this hypothesis excluded the most ``explicit'' of discrete subgroups of $PSL(2,\R )$, the modular
group $\Gamma =PSL(2,\Z )$.  The definition provided in this paper
is clearly available in this case, and we devote the rest of this section to its consideration.

Two elements $r,s\in \R\cup \{\infty\}\approx \SI^{1}$
are called {\em equivalent} if there exists $A\in PSL(2,\Z )$ such that
$A(r)=s$.  Every equivalence class $[r]$ of extended reals
corresponds to a leaf
$L_{[r]}$ of $\mathcal{F}$, and since all $PSL(2,\Z )$-orbits in $\SI^{1}$ are dense, 
all leaves are dense.  If $L_{[r]}$ is
isomorphic to the punctured hyperbolic disk $\D^{\ast}$, then $[r]$ is
quadratic over $\Q$.  Otherwise, $L_{[r]}$ is isomorphic to $\HP^{2}$.

Let us consider the leaf $L=L_{[0]}\cong\D^{\ast}$ covered by $\HP\times \{
0\}$.  Choose $x\in L$ and a transversal $T$ through $x$ that is a
fiber with respect to the projection onto the modular surface $\Sigma$.  We assume
that the lift $\tilde{x}$ of $x$ to $\HP^{2}$ is not an elliptic point
for the action of $PSL(2,\Z )$.
Define ${\sf f}$ to be the lift of a frame on $\Sigma$ based at the projection of $T$.  
As before, we denote by $\tilde{\sf f}$ the lift of ${\sf f}$ to $\widetilde{T}_{0}\subset\HP^{2}$
and by $\tilde{\sf f}_{\tilde{y}}$ its value at $\tilde{y}\in \widetilde{T}_{0}$.
Note that for $A\in  PSL(2,\R )$, $A_{\ast} \tilde{\sf f}_{\tilde{x}}=A_{\ast} \tilde{\sf f}_{\tilde{y}}$
if and only if $A\in \Gamma=PSL(2,\Z )$.
A sequence $\{ A_{\alpha}\}\subset PSL(2,\Z )$ defines an ${\sf f}$-diophantine
approximation $\Leftrightarrow$ $(A_{\alpha}\tilde{x},0)$ projects to a sequence in $T$ converging
to $x$ $\Leftrightarrow$ $(\tilde{x},A_{\alpha}^{-1}(0))$ projects to  a sequence in $T$ converging
to $x$ $\Leftrightarrow$ $A_{\alpha}^{-1}(0)\rightarrow 0$ in $\SI^{1}$.
Note that for any sequence $\{ \gamma_{n_{\alpha}}\}$ in the deck group of
$\HP^{2}\rightarrow L$,
\[ \left\{\gamma_{n}\;=\; \left( \left.\begin{array}{cc}
                              1 & 0 \\
                              n & 1
                              \end{array}\right)\;\right|\;\;n\in\Z     \right\} ,\]
 the sequence $\{\gamma_{\alpha}A_{\alpha}\}$ also defines an ${\sf f}$-diophantine
 approximation.                             
The fundamental germ $\fgrm (\mathcal{F},x,{\sf f})$ is then formed from the associated
sequences $\{ A_{\alpha}B_{\alpha}^{-1}\}$ where $\{B_{\alpha}\}$ is another
${\sf f}$-diophantine approximation.

\begin{note}  The sequences
$\{ b_{\alpha}/d_{\alpha}=A^{-1}_{\alpha}(0)\}$ are 
{\em hyperbolic diophantine approximations}, as defined for example in \cite{Be-Do}.  
See  \cite{Ge1} for more on this point.   In the case at hand, they give bad diophantine 
approximations
of $0$
whenever $b_{\alpha}\rightarrow\infty$, in the sense that it is never true that for some $c>0$
and almost all $\alpha$,
\[  \left| 0-\frac{b_{\alpha}}{d_{\alpha}}\right|\;\;
<\;\;\frac{c}{d^{2}_{\alpha}}.  \]
\end{note}

The ${\sf f}$-diophantine approximations are not stable with respect
to the operation of inversion.  Indeed, let $r\in\R$ be any real 
number, $\{ m_{\alpha}/n_{\alpha}\}$ a sequence of rationals 
(written in 
lowest terms)
converging to $r$.  Let $M_{\alpha}, N_{\alpha}$  be such that
$m_{\alpha}M_{\alpha}-n_{\alpha}N_{\alpha}=1$.   Assume that the indexing
is such that $\alpha+ N_{\alpha}/m_{\alpha}\rightarrow\infty$ as $\alpha\rightarrow\infty$.
Then the sequence $\{ 
X_{\alpha}\}$, 
\begin{equation}\label{unstable}   X_{\alpha}\;\;=\;\; \left(\begin{array}{cc}
                                                     -(\alpha m_{\alpha}+N_{\alpha} )&  
                                                    -m_{\alpha}  \\
                                                    \alpha 
                                                   n_{\alpha}+M_{\alpha}  
                                                   &   n_{\alpha} 
						   \end{array}
\right) \;\;\in\;\; PSL(2, \Z ) ,
\end{equation}
satisfies $X_{\alpha}^{-1}(0)\rightarrow 0$,  but 
$X_{\alpha}(0)\rightarrow -r$.  Using this fact, we can now show

\begin{theo} As a set, 
    \[ \fgrm (\mathcal{F}, x, {\sf f})\;\; =\;\; PSL(2,\bast \Z ).\]
    \end{theo}
\begin{proof}  Let $\{ A_{\alpha}\}$ be any sequence in $PSL 
(2,\bast\Z )$.   Then after passing to a subsequence if necessary, we
find $A^{-1}_{\alpha}(0)\rightarrow r$ for some $r\in \R\cup 
\{\infty\}$.  Note that $r$ is 
independent of the class of $\{ A_{\alpha}\}$ in $PSL(2,\bast \Z )$.
We may choose $\{ n_{\alpha}\}\subset\Z$ so that 
$\gamma_{n_{\alpha}}A^{-1}_{\alpha}(0)\rightarrow 0$.  Hence
\[ \{ A_{\alpha}\}\;\; =\;\;\{ A_{\alpha}\gamma_{-n_{\alpha}}\}
\cdot \{\gamma^{-1}_{-n_{\alpha}}\}\]
defines an element of $\fgrm (\mathcal{F}, x, {\sf f})$.
\end{proof}

It is not difficult to see that with respect to its action on the germ universal
cover $\gcov{\mathcal{F}}$,
$\fgrm (\mathcal{F},x,{\sf f})$  is not a group.  Indeed, the class of the sequence
$\{ X_{\alpha}^{-1}\}$, where  $\{ X_{\alpha}\}$
is the sequence appearing in (\ref{unstable}), is not defined on $\widetilde{L}$.

\section{Mapping Class Group of the Algebraic Universal Cover of a Surface}

In this section, we use the fundamental germ to
prove a Nielsen type theorem for the algebraic universal cover of a closed
surface.  We begin by recalling a few facts, referring the reader to \cite{Ge2} for details.

Let $\Sigma$ be a closed surface and let $\mathcal{G}=\{ 
G_{\alpha}\}$ be
the set of all normal finite index subgroups.  For each $G_{\alpha}$,
there exists a covering $\sigma_{\alpha}:\Sigma_{\alpha}\rightarrow \Sigma$ defined
by the condition that $\pi_{1}\Sigma_{\alpha}$ maps isomorphically
onto $G_{\alpha}$.  If $G_{\alpha}\subset G_{\beta}$, there is a unique
covering 
$s_{\alpha\beta}:\Sigma_{\alpha}\rightarrow\Sigma_{\beta}$ 
for which $\sigma_{\alpha}=\sigma_{\beta}\circ s_{\alpha\beta}$.
Hence the collection of $\sigma_{\alpha}$ and $s_{\alpha\beta}$ forms
an inverse system of surfaces by covering maps.

\begin{defi} The {\bf algebraic universal cover} of $\Sigma$ is the 
inverse limit 
\[    \widehat{\Sigma}\;\; =\;\;\lim_{\longleftarrow} \Sigma_{\alpha} .     \]
\end{defi}  

If $\sigma :Z\rightarrow \Sigma$ is any finite covering, 
then $\sigma$ lifts to a homeomorphism
\[ \hat{\sigma}:\widehat{Z}\rightarrow\widehat{\Sigma}.\]
Thus the algebraic universal cover depends only 
on the type of $\Sigma$ (elliptic, parabolic, hyperbolic).
In fact, there are only two non-trivial examples of algebraic universal
covers of closed surfaces: that
of the torus and that
of a surface of hyperbolic type.

The inverse limit \[ \hat{\pi}_{1}\Sigma\;\;=\;\; 
\lim_{\longleftarrow}\big(\pi_{1}\Sigma\big)/G_{\alpha}\] is 
a Cantor group called the {\em profinite completion} of $\pi_{1}\Sigma$.
The homomorphism 
$i :\pi_{1}\Sigma\rightarrow
\hat{\pi}_{1}\Sigma$ induced by the system of projections 
$\pi_{1}\Sigma\rightarrow \pi_{1}\Sigma/G_{\alpha}$ has dense image.  Define
a representation 
\begin{eqnarray*}  & \varsigma :\pi_{1}\Sigma\longrightarrow 
{\sf Homeo}(\hat{\pi}_{1}\Sigma ) & \\
                   & \varsigma_{\gamma}(\hat{g} ) = \hat{g}\cdot i(\gamma)^{-1} &     
\end{eqnarray*}  
for $\gamma\in\pi_{1}\Sigma$ and $\hat{g}\in\hat{\pi}_{1}\Sigma$.
Then 
we may 
identify
$\widehat{\Sigma}$ with the suspension of $\varsigma$:
\[ \widehat{\Sigma}\;\;\approx\;\;    
\big( \widetilde{\Sigma}\times\hat{\pi}_{1}\Sigma\big)\big/\pi_{1}\Sigma .\]
With this identification, we see that $\widehat{\Sigma}$ is a surface lamination with
Cantor transversals homeomorphic to $\hat{\pi}_{1}\Sigma$, that is, a solenoid.  Moreover, it can also be seen from this presentation
that every leaf $L$ of $\widehat{\Sigma}$ satisfies 
\[ \pi_{1}L\;\; \cong\;\; \bigcap G_{\alpha}.  \]
However for closed surfaces, $\bigcap G_{\alpha}=1$, so here, $L$ is simply connected.
Each leaf $L$ is dense and a path-component of $\widehat{\Sigma}$.
For every $\alpha$, the pre-image of the projection map 
$\widehat{\Sigma}\rightarrow \Sigma_{\alpha}$
is a fiber transversal, homeomorphic to 
$\hat{\pi}_{1}\Sigma_{\alpha}\cong\widehat{G}_{\alpha}$.

Now let $L$ be a fixed leaf of $\widehat{\Sigma}$.
\begin{defi}  The {\bf leafed mapping class group} of $\widehat{\Sigma}$ is
\[ {\sf MCG}(\widehat{\Sigma},L)\;\; =\;\; {\sf Homeo}(\widehat{\Sigma},L)/\simeq ,
\]
where $\simeq$ is the relation of homotopy of homeomorphisms.  
\end{defi}
We denote by $[h]$ the mapping class
associated to a homeomorphism $h$.

Let $G$ be a group. 

\begin{defi} The  {\bf virtual automorphism group} of $G$ is 
\[ {\sf Vaut}(G)\;\;=\;\;
\Big\{ \phi: H\rightarrow H'\; \Big|\;\;
\phi\text{ an 
isomorphism and }H,\, H'\text{ finite index subgroups of } G
  \Big\}\Big/\sim , 
\]
where $\phi_{1}\sim\phi_{2}$ if there exists 
$H''<G$ of finite index, contained in ${\sf 
Dom}(\phi_{1})\cap {\sf Dom}(\phi_{2})$ and such that $\phi_{1}|_{H''}=
\phi_{2}|_{H''}$. 
\end{defi}

Note that the equivalence relation $\sim$ is precisely what is needed
to make composition of virtual automorphisms well-defined.
We point out also that if $H<G$ is of finite index, then ${\sf Vaut}(H)\cong {\sf Vaut}(G)$.

\begin{theo}\label{Nielsentheorem} 
${\sf MCG}(\widehat{\Sigma},L)\;\;\cong\;\;{\sf Vaut}(\pi_{1}\Sigma )$.
\end{theo}

\begin{note}  We first learned the statement of this theorem in a conversation
with D. Sullivan in 1995.  The first proof appeared in the 
thesis of C. Odden \cite{Od}.
\end{note}

\begin{proof} Define a homomorphism
\[\Theta :{\sf Vaut}(\pi_{1}\Sigma)\longrightarrow {\sf MCG}(\widehat{\Sigma},L)\]
as follows.  Given $\phi : G_{1}\rightarrow G_{2}$ an isomorphism
of finite index subgroups of $\pi_{1}\Sigma$, we may find
covers $\sigma_{1}, \sigma_{2}:\Sigma'\rightarrow \Sigma$ -- indexed 
by $G_{1}$ and $G_{2}$ -- so 
that 
\[     (\sigma_{2})_{\ast} \circ (\sigma_{1})_{\ast}^{-1} \;\;=\;\; \phi  :
\] this follows from the classical Nielsen theorem.
Then we define
\[   \Theta (\phi ) \;\; =\;\; 
[\hat{\sigma}_{2}\circ\hat{\sigma}_{1}^{-1}]\]
where for $i=1,2$,  $\hat{\sigma}_{i}:\widehat{\Sigma}'\rightarrow \widehat{\Sigma}$ is
the algebraic universal cover lift.  If $G'< {\sf Dom}(\phi )$, then 
$\Theta (\phi|_{G'})\;\; =\;\; \Theta (\phi )$,
since $\Theta (\phi|_{G'})$ is defined by the pair $\sigma_{i}\circ\sigma$, 
$i=1,2$, where
$\sigma:\Sigma''\rightarrow \Sigma'$ is a cover for which $\sigma_{1}\circ\sigma$ is indexed by $G'$.
Thus $\Theta$ is a well-defined homomorphism.

\begin{clai} $\Theta$ {\em is onto.}
\end{clai}

Let $h:(\widehat{\Sigma},L)\rightarrow (\widehat{\Sigma},L)$ 
be a homeomorphism.  After performing
an isotopy, we may arrange that 
$h$ fixes a point $x$ and fiber transversal $T$ containing $x$.  
Without loss of generality, we may assume that $T$
is a fiber transversal over $\Sigma$.  Due to the suspension structure,
$T\approx\hat{\pi}\Sigma$: fix this identification 
so that $x\mapsto 1$ and
$L\cap T\mapsto\pi_{1}\Sigma$.  
Since $h(L)=L$, we obtain a bijection
\[ h_{\ast}:\pi_{1}\Sigma\longrightarrow \pi_{1}\Sigma  \]
in which $h(1)=1$.

Suppose that for each $G_{\alpha}<\pi_{1}\Sigma$, $h_{\ast}|G_{\alpha}$ is not homomorphic.
This means that for every $\alpha$,
there exists $\gamma_{\alpha},\gamma_{\alpha}'\in G_{\alpha}$ so that
\begin{equation}\label{nothomomorphic}  
h_{\ast}\big(\gamma_{\alpha}\cdot\gamma_{\alpha}'\big)  \not= 
h_{\ast}\big(\gamma_{\alpha}\big)\cdot h_{\ast}\big(\gamma_{\alpha}'\big).       
\end{equation}
Assuming that $\widehat{\Sigma}$ has been equipped with a 
hyperbolic metric, say lifted from $\Sigma$, then
the sequences $\{ \gamma_{\alpha}\}$, $\{\gamma_{\alpha}'\}$ define elements of the 
fundamental germ \[ \bast\gamma,\; \bast\gamma'\;\;\in\;\;
\fgrm (\widehat{\Sigma},x,{\sf f})\]
where ${\sf f}$ is a frame field lifted from a frame on $\Sigma$.  But this
fundamental germ is a subgroup of the mother germ $\fgrm (\widehat{\Sigma})$.
By Theorem~\ref{functor}, $h$ induces a groupoid isomorphism
\[  \gzg{h}_{\ast}:\fgrm(\widehat{\Sigma})\;\longrightarrow\; 
\fgrm(\widehat{\Sigma}),\] 
and so we must have  
\[  \gzg{h}_{\ast} \big(\bast\gamma\cdot\bast\gamma'\big)\;\; =\;\; 
\gzg{h}_{\ast}\big(\bast\gamma\big)\cdot
\gzg{h}_{\ast}\big(\bast\gamma'\big).\]
This contradicts equation (\ref{nothomomorphic}).  
Thus $h_{\ast}$ defines an isomorphism when restricted to some $G_{\alpha}$,
and this isomorphism determines an element $\phi\in {\sf Vaut}(\pi_{1}\Sigma )$.
Note that that $\phi$ does not depend on the isotopy used to ensure $h(x)=x$
since the holonomy group of $\widehat{\Sigma}$ at any point is trivial.
Choose $\sigma_{i}:\Sigma'\rightarrow\Sigma$, $i=1,2$, so that
$\Theta (\phi )=[\hat{\sigma}_{2}\circ\hat{\sigma}_{1}^{-1}]$.
To simplify notation, we write $h_{0}=\hat{\sigma}_{2}\circ\hat{\sigma}_{1}^{-1}$.

Recall that since $\widehat{\Sigma}$ is compact with hyperbolic
leaves, the germ universal cover of $\widehat{\Sigma}$ is $\bstar\HP^{2}$.
The homeomorphisms
$h$ and $h_{0}$ lift to the standard bijections
$\bstar h$ and $\bstar h_{0}$ of $\bstar\HP^{2}$
sharing the same equivariance with respect to the action of the mother germ
$\fgrm (\widehat{\Sigma})$.  
In particular, they act identically on the set
of galaxies of $\bstar\HP^{2}$.  For this reason, we may choose a germ universal
cover topology for $\bstar\HP^{2}$ with respect to which both 
$\bstar h$ and $\bstar h_{0}$ are homeomorphisms.

Define a homotopy $\bstar H_{t}$ from $\bstar h$ to 
$\bstar h_{0}$ as follows.  For each
$\bstar z\in\bstar \HP^{2}$, 
$\bstar H_{t}\big(\bstar z\big)$ is the point 
subdividing
the hyperbolic geodesic connecting $\bstar h\big(\bstar z\big)$ to
$\bstar h_{0}\big(\bstar z\big)$ into the proportion $t:1-t$.  
By construction,
$\bstar H_{t}$ has the same equivariance as $\bstar h$ and 
$\bstar h_{0}$ and is in particular continuous.  
Since its initial and final maps are standard,
so is $\bstar H_{t}$.  By Proposition~\ref{descend}, it 
descends to a homotopy $H_{t}$ of $h$ and $h_{0}$. 
It follows that $[h]=[h_{0}]=\Theta (\phi )$,
and $\Theta$ is 
onto.

\begin{clai}  $\Theta$ is one-to-one.
\end{clai}

If not, then there exists $\phi\not=$ the identity map with $\Theta (\phi )=1$.
But then $\Theta (\phi )$ would have to induce the identity map on the
mother germ;
by construction, this can only happen if $\phi$ is trivial.
\end{proof}

Let ${\sf Mod}(\widehat{\Sigma},L)$ be the {\em Teichm\"{u}ller modular
group} of the pair $(\widehat{\Sigma},L)$: the group of homotopy
classes of quasiconformal homeomorphisms of $\widehat{\Sigma}$ that
preserve $L$.

\begin{coro}\label{Teichmapgroup}  ${\sf Mod}(\widehat{\Sigma},L)={\sf MCG}(\widehat{\Sigma},L)$.
\end{coro}

\begin{proof}  This follows from the proof of Theorem~\ref{Nielsentheorem}
and the fact that every finite cover
of compact Riemann surfaces is homotopic to a quasiconformal cover.
\end{proof}

Theorem~\ref{Nielsentheorem} can be used to formulate
the following conjectural Nielsen-type theorem.  
Given $x\in L$, the fundamental germ $\fgrm (\widehat{\Sigma},x)$ is made
up of all sequences $\{ \gamma_{\alpha}\}$ converging with respect to
the lattice of finite index normal subgroups $G$ of $\pi_{1}\Sigma$, so
\[ \fgrm (\widehat{\Sigma},x)\;\;\cong \;\;
\bigcap_{[\pi_{1}\Sigma :G ]<\infty}\bast G\;\;\subset\;\; \bast\pi_{1}\Sigma .\]
It follows then that there is a monomorphism 
${\sf Vaut}(\pi_{1})\hookrightarrow {\sf Aut}(\fgrm (\widehat{\Sigma},x))$,
which descends to 
${\sf Vaut}(\pi_{1})\hookrightarrow {\sf Out}(\fgrm (\widehat{\Sigma},x))$
upon passage to the quotient.  This latter map is also a monomorphism,
since no nontrivial virtual automorphism $\phi$ 
can induce on $\fgrm (\widehat{\Sigma},x)$
an inner automorphism.  For otherwise, $\phi$ would have to be inner on
some subgroup $H$, hence all, which is only possible if $\phi$ is trivial.
In view of these remarks we 

\begin{conj}  The monomorphism ${\sf MCG}(\widehat{\Sigma},L)\hookrightarrow
{\sf Out}(\fgrm (\widehat{\Sigma},x))$ is an isomorphism.
\end{conj}

We end this section by explaining
the importance of Theorem~\ref{Nielsentheorem} and Corollary~\ref{Teichmapgroup} in giving a genus
independent reformulation the Ehrenpreis conjecture. 
The classical Ehrenpreis conjecture is:

\begin{quote}
\textit{Given two closed hyperbolic surfaces $\Sigma_{1}$ and $\Sigma_{2}$
and $\epsilon>0$,
there exist finite, locally isometric covering surfaces $Z_{1}$ and $Z_{2}$
of each which are $(1+\epsilon )$-quasiisometric.}
\end{quote}

We then have the following
equivalent,
genus independent version:

\begin{quote} \textit{Every orbit of the action of
${\sf Mod}(\widehat{\Sigma},L)$ on $\mathcal{T}(\widehat{\Sigma})$ 
is dense.}
\end{quote}

In other words, the genus independent version says
that, although the moduli space
\[ \mathcal{T}(\widehat{\Sigma})/{\sf Mod}(\widehat{\Sigma},L)\]
is uncountable, it has the ``topology of a point'' ({\em i.e.}\ the coarse
topology).  If affirmed, the Ehrenpreis
conjecture would thus provide an explanation for the jump between the existence
of moduli (dimension 2) and rigidity (dimension 3 and higher) in hyperbolic
geometry.  
See the articles \cite{BNS}, \cite{Ge2} for more discussion.

\bibliographystyle{amsalpha}

\end{document}